\documentclass[1pt]{amsart}

\usepackage{amsmath,amsthm, amscd, amssymb, amsfonts}

\usepackage{hhline}

\newcommand{\Ind}{\operatorname{Ind}}
\newcommand{\ord}{\operatorname{ord}}
\newcommand{\oper}{\#}

\DeclareMathOperator*{\prode}{\boxtimes}

\newcommand\toba{{\mathfrak B }}
\newcommand\trasp{\pi}
\newcommand{\gr}{\operatorname{gr}}
\newcommand{\trid}{\triangleright}

\newcommand\compvert[2]{\genfrac{}{}{-1pt}{0}{#1}{#2}}

\newcommand\mvert[2]{\begin{tiny}\begin{matrix}#1\vspace{-4pt}\\#2\end{matrix}\end{tiny}}
\DeclareMathOperator*{\Tim}{\times}
\newcommand{\Times}[2]{\sideset{_#1}{_#2}\Tim}

\newcommand{\cx}{{\daleth}}

\newcommand{\lu}{{\nu}}
\newcommand{\ld}{{\ell}}

\newcommand{\abofe}{{A(\lu, \ld)}}

\newcommand{\bbofe}{{B(\lu, \ld)}}

\newcommand{\proy}{\Pi}

\newcommand{\Lc}{{\mathcal L}}
\newcommand{\Y}{{\mathcal Y}}
\newcommand{\R}{{\mathcal R}}
\newcommand{\W}{{\mathcal W}}
\newcommand{\Zc}{{\mathcal Z}}
\newcommand{\daga}{{\dagger}}

\newcommand{\acts}{{\rightharpoondown}}

\newcommand{\ku}{\mathbb C}
\newcommand{\K}{{\mathcal K}}
\newcommand{\Z}{{\mathbb Z}}
\newcommand{\N}{{\mathbb N}}
\newcommand{\I}{{\mathbb I}}

\newcommand{\G}{{\mathbb G}}
\newcommand{\M}{{\mathcal M}}
\newcommand{\Q}{{\mathcal Q}}
\newcommand{\F}{{\mathcal F}}
\newcommand{\C}{{\mathcal C}}
\newcommand{\D}{{\mathcal D}}
\newcommand{\Ec}{{\mathbf E}}
\newcommand{\Ee}{{\mathcal E}}
\newcommand{\Kc}{{\mathbf K}}

\newcommand{\uno}{{\bf 1}}
\newcommand{\uv}{{\bf 1_v}}
\newcommand{\uh}{{\bf 1_h}}
\newcommand{\m}{\mathcal{M}}
\newcommand{\n}{\mathcal{N}}

\newcommand{\Dc}{{\mathbf D}}
\newcommand{\B}{{\mathcal B}}
\newcommand{\T}{{\mathcal T}}
\newcommand{\Hc}{{\mathcal H}}
\newcommand{\Vc}{{\mathcal V}}
\newcommand{\Pc}{{\mathcal P}}
\newcommand{\Oc}{{\mathcal O}}
\newcommand{\ydh}{{}^H_H\mathcal{YD}}

\newcommand{\Ss}{{\mathcal S}}
\newcommand{\Vect}{\operatorname{Vec}}
\newcommand{\End}{\operatorname{End}}
\newcommand{\Aut}{\operatorname{Aut}}
\newcommand{\Int}{\operatorname{Int}}
\newcommand{\Ext}{\operatorname{Ext}}
\newcommand\tr{\operatorname{tr}}
\newcommand\Rep{\operatorname{Rep}}

\newcommand\card{\operatorname{card}}
\newcommand{\unosigma}{{\bf 1}^{\sigma}}
\newcommand\tot{\operatorname{tot}}
\newcommand\sgn{\operatorname{sgn}}
\newcommand\ad{\operatorname{ad}}
\newcommand\Hom{\operatorname{Hom}}
\newcommand\opext{\operatorname{Opext}}
\newcommand\Tot{\operatorname{Tot}}
\newcommand\Map{\operatorname{Map}}
\newcommand{\fde}{{\triangleright}}
\newcommand{\fiz}{{\triangleleft}}
\newcommand{\wC}{\widehat{C}}
\newcommand{\la}{\langle}
\newcommand{\ra}{\rangle}
\newcommand{\Comod}{\mbox{\rm Comod\,}}
\newcommand{\Mod}{\mbox{\rm Mod\,}}
\newcommand{\Sg}{{\mathfrak S}}
\newcommand{\funcion}{\mathfrak p}
\newcommand{\tauo}{{\widehat\tau}}
\newcommand{\prin}{t}
\newcommand{\fin}{b}
\newcommand{\pri}{r}
\newcommand{\fine}{l}
\newcommand\rh{\sim_{H}}
\newcommand\rv{\sim_{V}}
\newcommand\rd{\sim_{D}}

\newcommand\V{\operatorname{Vec}}

\theoremstyle{plain}

\renewcommand{\baselinestretch}{1.2}

\title{The character tables of  centralizers in  Weyl Groups of $E_8$ I}
\author{ \small Shouchuan Zhang,    \ \ Peng Wang, \ \ Jing Cheng,\ \ Hui Yang
}
\address{ Mathematics Department of Hunan University,\newline \indent Changsha China,
410082, E-mail: z9491@yahoo.com.cn }
\date{}

\begin{document}
\newtheorem{Proposition}{\quad Proposition}[section]
\newtheorem{Theorem}{\quad Theorem}
\newtheorem{Definition}[Proposition]{\quad Definition}
\newtheorem{Corollary}[Proposition]{\quad Corollary}
\newtheorem{Lemma}[Proposition]{\quad Lemma}
\newtheorem{Example}[Proposition]{\quad Example}

\maketitle \addtocounter{section}{-1}

\numberwithin{equation}{section}

\date{}


\begin {abstract} To classify the finite dimensional pointed Hopf
algebras with Weyl group $G$ of $E_8$, we obtain  the
representatives of conjugacy classes of $G$ and  all character
tables of centralizers of these representatives by means of software
GAP. In this paper we only list character table 1-- 28.

\vskip0.1cm 2000 Mathematics Subject Classification: 16W30, 68M07

keywords: GAP, Hopf algebra, Weyl group, character.
\end {abstract}

\section{Introduction}\label {s0}

This article is to contribute to the classification of
finite-dimensional complex pointed Hopf algebras  with Weyl groups
of $E_8$.
 Many papers are about the classification of finite dimensional
pointed Hopf algebras, for example,  \cite{AS98b, AS02, AS00, AS05,
He06, AHS08, AG03, AFZ08, AZ07,Gr00,  Fa07, AF06, AF07, ZZC04, ZCZ08}.

In these research  ones need  the centralizers and character tables
of groups. In this paper we obtain   the representatives of
conjugacy classes of Weyl groups of $E_8$  and all
character tables of centralizers of these representatives by means
of software GAP. In this paper we only list character table 1-- 28.

By the Cartan-Killing classification of simple Lie algebras over complex field the Weyl groups
to
be considered are $W(A_l), (l\geq 1); $
 $W(B_l), (l \geq 2);$
 $W(C_l), (l \geq 2); $
 $W(D_l), (l\geq 4); $
 $W(E_l), (8 \geq l \geq 6); $
 $W(F_l), ( l=4 );$
$W(G_l), (l=2)$.

It is otherwise desirable to do this in view of the importance of Weyl groups in the theories of
Lie groups, Lie algebras and algebraic groups. For example, the irreducible representations of
Weyl groups were obtained by Frobenius, Schur and Young. The conjugace classes of  $W(F_4)$
were obtained by Wall
 \cite {Wa63} and its character tables were obtained by Kondo \cite {Ko65}.
The conjugace classes  and  character tables of $W(E_6),$ $W(E_7)$ and $W(E_8)$ were obtained
by  Frame \cite {Fr51}.
Carter gave a unified description of the conjugace classes of
 Weyl groups of simple Lie algebras  \cite  {Ca72}.

\section {Program}

 By using the following
program in GAP,  we  obtain  the representatives of conjugacy
classes of Weyl groups of $E_6$  and all character tables of
centralizers of these representatives.

gap$>$ L:=SimpleLieAlgebra("E",6,Rationals);;

gap$>$ R:=RootSystem(L);;

 gap$>$ W:=WeylGroup(R);Display(Order(W));

gap $>$ ccl:=ConjugacyClasses(W);;

gap$>$ q:=NrConjugacyClasses(W);; Display (q);

gap$>$ for i in [1..q] do

$>$ r:=Order(Representative(ccl[i]));Display(r);;

$>$ od; gap

$>$ s1:=Representative(ccl[1]);;cen1:=Centralizer(W,s1);;

gap$>$ cl1:=ConjugacyClasses(cen1);

gap$>$ s1:=Representative(ccl[2]);;cen1:=Centralizer(W,s1);;

 gap$>$
cl2:=ConjugacyClasses(cen1);

gap$>$ s1:=Representative(ccl[3]);;cen1:=Centralizer(W,s1);;

gap$>$ cl3:=ConjugacyClasses(cen1);

 gap$>$
s1:=Representative(ccl[4]);;cen1:=Centralizer(W,s1);;

 gap$>$
cl4:=ConjugacyClasses(cen1);

 gap$>$
s1:=Representative(ccl[5]);;cen1:=Centralizer(W,s1);;

gap$>$ cl5:=ConjugacyClasses(cen1);

 gap$>$
s1:=Representative(ccl[6]);;cen1:=Centralizer(W,s1);;

gap$>$ cl6:=ConjugacyClasses(cen1);

gap$>$ s1:=Representative(ccl[7]);;cen1:=Centralizer(W,s1);;

gap$>$ cl7:=ConjugacyClasses(cen1);

gap$>$ s1:=Representative(ccl[8]);;cen1:=Centralizer(W,s1);;

gap$>$ cl8:=ConjugacyClasses(cen1);

gap$>$ s1:=Representative(ccl[9]);;cen1:=Centralizer(W,s1);;

gap$>$ cl9:=ConjugacyClasses(cen1);

gap$>$ s1:=Representative(ccl[10]);;cen1:=Centralizer(W,s1);;

 gap$>$
cl10:=ConjugacyClasses(cen1);

gap$>$ s1:=Representative(ccl[11]);;cen1:=Centralizer(W,s1);;

gap$>$ cl11:=ConjugacyClasses(cen1);

gap$>$ s1:=Representative(ccl[12]);;cen1:=Centralizer(W,s1);;

gap$>$ cl2:=ConjugacyClasses(cen1);

gap$>$ s1:=Representative(ccl[13]);;cen1:=Centralizer(W,s1);;

 gap$>$
cl13:=ConjugacyClasses(cen1);

 gap$>$
s1:=Representative(ccl[14]);;cen1:=Centralizer(W,s1);;

gap$>$ cl14:=ConjugacyClasses(cen1);

gap$>$ s1:=Representative(ccl[15]);;cen1:=Centralizer(W,s1);;

gap$>$ cl15:=ConjugacyClasses(cen1);

gap$>$ s1:=Representative(ccl[16]);;cen1:=Centralizer(W,s1);;

gap$>$ cl16:=ConjugacyClasses(cen1);

gap$>$ s1:=Representative(ccl[17]);;cen1:=Centralizer(W,s1);;

gap$>$ cl17:=ConjugacyClasses(cen1);

gap$>$ s1:=Representative(ccl[18]);;cen1:=Centralizer(W,s1);;

gap$>$ cl18:=ConjugacyClasses(cen1);

gap$>$ s1:=Representative(ccl[19]);;cen1:=Centralizer(W,s1);;

gap$>$ cl19:=ConjugacyClasses(cen1);

gap$>$ s1:=Representative(ccl[20]);;cen1:=Centralizer(W,s1);;

gap$>$ cl20:=ConjugacyClasses(cen1);

gap$>$ s1:=Representative(ccl[21]);;cen1:=Centralizer(W,s1);;

 gap$>$
cl21:=ConjugacyClasses(cen1);

gap$>$ s1:=Representative(ccl[22]);;cen1:=Centralizer(W,s1);;

 gap$>$
cl22:=ConjugacyClasses(cen1);

gap$>$ s1:=Representative(ccl[23]);;cen1:=Centralizer(W,s1);;

gap$>$ cl23:=ConjugacyClasses(cen1);

gap$>$ s1:=Representative(ccl[24]);;cen1:=Centralizer(W,s1);;

$>$ cl24:=ConjugacyClasses(cen1);

gap$>$ s1:=Representative(ccl[25]);;cen1:=Centralizer(W,s1);

gap$>$ cl25:=ConjugacyClasses(cen1);

gap$>$ for i in [1..q] do

$>$ s:=Representative(ccl[i]);;cen:=Centralizer(W,s);;

$>$ char:=CharacterTable(cen);;Display (cen);Display(char);

 $>$ od;

The programs for Weyl groups of $E_7$, $E_8$, $F_4$ and $ G_2$ are
similar.

\section {$E_8$}

 The generators of $G$ are:\\\\
 $\left(
\begin{array}{cccccccc} -1&0&1&0&0&0&0&0
\\ 0&1&0&0&0&0&0&0 \\ 0&0&1&0&0&0&0&0 \\ 0&0&0&
1&0&0&0&0 \\ 0&0&0&0&1&0&0&0 \\ 0&0&0&0&0&1&0&0 \\
  0&0&0&0&0&0&1&0 \\ 0&0&0&0&0&0&0&1
\end{array}\right)$,
  $\left(\begin{array}{cccccccc} 1&0&0&0&0&0&0&0 \\ 0&-1&0&1&0&0&0&0 \\ 0&0&1&0&0&0&0&0 \\
        0&0&0&1&0&0&0&0 \\ 0&0&0&0&1&0&0&0 \\ 0&0&0&0&0&1&0&0 \\
        0&0&0&0&0&0&1&0 \\ 0&0&0&0&0&0&0&1 \end{array}\right)$,
  $\left(\begin{array}{cccccccc} 1&0&0&0&0&0&0&0 \\ 0&1&0&0&0&0&0&0 \\ 1&0&-1&1&0&0&0&0 \\
        0&0&0&1&0&0&0&0 \\ 0&0&0&0&1&0&0&0 \\ 0&0&0&0&0&1&0&0 \\
        0&0&0&0&0&0&1&0 \\ 0&0&0&0&0&0&0&1 \end{array}\right)$,\\
  $\left(\begin{array}{cccccccc} 1&0&0&0&0&0&0&0 \\ 0&1&0&0&0&0&0&0 \\ 0&0&1&0&0&0&0&0 \\
        0&1&1&-1&1&0&0&0 \\ 0&0&0&0&1&0&0&0 \\ 0&0&0&0&0&1&0&0 \\
        0&0&0&0&0&0&1&0 \\ 0&0&0&0&0&0&0&1 \end{array}\right)$,
  $\left(\begin{array}{cccccccc} 1&0&0&0&0&0&0&0 \\ 0&1&0&0&0&0&0&0 \\ 0&0&1&0&0&0&0&0 \\
        0&0&0&1&0&0&0&0 \\ 0&0&0&1&-1&1&0&0 \\ 0&0&0&0&0&1&0&0 \\
        0&0&0&0&0&0&1&0 \\ 0&0&0&0&0&0&0&1 \end{array}\right)$,
  $\left(\begin{array}{cccccccc} 1&0&0&0&0&0&0&0 \\ 0&1&0&0&0&0&0&0 \\ 0&0&1&0&0&0&0&0 \\
        0&0&0&1&0&0&0&0 \\ 0&0&0&0&1&0&0&0 \\ 0&0&0&0&1&-1&1&0 \\
        0&0&0&0&0&0&1&0 \\ 0&0&0&0&0&0&0&1 \end{array}\right)$,\\
  $\left(\begin{array}{cccccccc} 1&0&0&0&0&0&0&0 \\ 0&1&0&0&0&0&0&0 \\ 0&0&1&0&0&0&0&0 \\
        0&0&0&1&0&0&0&0 \\ 0&0&0&0&1&0&0&0 \\ 0&0&0&0&0&1&0&0 \\
        0&0&0&0&0&1&-1&1 \\ 0&0&0&0&0&0&0&1 \end{array}\right)$,
  $\left(\begin{array}{cccccccc} 1&0&0&0&0&0&0&0 \\ 0&1&0&0&0&0&0&0 \\ 0&0&1&0&0&0&0&0 \\
        0&0&0&1&0&0&0&0 \\ 0&0&0&0&1&0&0&0 \\ 0&0&0&0&0&1&0&0 \\
        0&0&0&0&0&0&1&0 \\ 0&0&0&0&0&0&1&-1 \end{array}\right).$

 The representatives of conjugacy classes of $G$ are:\\\\
$s_{1} := \left(\begin{array}{cccccccc} 1&0&0&0&0&0&0&0 \\ 0&1&0&0&
0&0&0&0 \\ 0&0&1&0&0&0&0&0 \\
        0&0&0&1&0&0&0&0 \\ 0&0&0&0&1&0&0&0 \\ 0&0&0&0&0&1&0&0 \\
        0&0&0&0&0&0&1&0 \\ 0&0&0&0&0&0&0&1 \end{array}\right)$,
  $s_{2} := \left(\begin{array}   {cccccccc} 1&1&1&-2&0&0&1&0 \\ 1&1&1&-2&0&-1&2&0 \\ 2&1&2&-3&0&-1&2&0 \\
        2&2&3&-4&-1&-1&3&0 \\ 1&1&3&-3&-1&-1&3&0 \\ 1&1&2&-2&-1&-1&2&1 \\
        1&0&1&-1&-1&0&1&1 \\ 1&0&0&0&-1&0&1&0 \end{array}\right)$,
  $s_{3} := \left(\begin{array}   {cccccccc} 1&-1&-1&0&1&0&0&1 \\ 1&-2&-2&1&1&1&-1&1 \\ 2&-2&-2&0&2&1&-1&1 \\
        3&-3&-4&1&3&1&-2&2 \\ 3&-3&-3&1&2&1&-2&2 \\ 3&-2&-3&1&1&1&-1&1 \\
        2&-1&-2&1&0&1&-1&1 \\ 1&0&-1&0&0&1&-1&1
        \end{array}\right)$,\\
  $s_{4} := \left(\begin{array}   {cccccccc} 0&1&0&-1&0&0&0&2 \\ -1&1&1&-1&0&-1&0&3 \\ 0&2&1&-2&0&-1&0&4 \\
        -1&3&2&-3&0&-1&-1&6 \\ -1&2&2&-2&-1&0&-1&5 \\ -1&2&1&-1&-1&0&-1&4 \\ 0&1&0&0&-1&0&-1&3 \\
        0&1&0&0&-1&0&0&1 \end{array}\right)$,
  $s_{5} := \left(\begin{array}   {cccccccc} 0&1&-1&0&-1&1&1&-1 \\ 0&0&-1&0&-1&2&1&-2 \\
        1&1&-2&0&-2&3&1&-2 \\ 1&1&-2&0&-3&4&2&-4 \\ 1&0&-2&1&-3&3&2&-3 \\
        1&0&-2&1&-2&2&1&-2 \\ 0&0&-1&1&-2&2&0&-1 \\ 0&0&-1&1&-1&1&0&-1
        \end{array}\right)$,  \\
  $s_{6} := \left(\begin{array}   {cccccccc} 0&0&1&0&-2&1&1&0 \\ -1&1&2&0&-3&1&1&0 \\ -1&1&3&-1&-3&1&2&0 \\
        -2&2&4&-1&-5&2&3&-1 \\ -1&2&3&-1&-4&2&2&-1 \\
        -1&2&2&-1&-3&2&1&0 \\ 0&2&1&-1&-2&1&1&0 \\ 0&1&0&0&-1&0&1&0 \end{array}\right)$,
  $s_{7} := \left(\begin{array}   {cccccccc} -1&0&0&0&0&0&0&0 \\ 0&-1&0&0&0&0&0&0 \\ 0&0&-1&0&0&0&0&0 \\
        0&0&0&-1&0&0&0&0 \\ 0&0&0&0&-1&0&0&0 \\ 0&0&0&0&0&-1&0&0 \\
        0&0&0&0&0&0&-1&0 \\ 0&0&0&0&0&0&0&-1 \end{array}\right)$,
        \\
  $s_{8} := \left(\begin{array}   {cccccccc} 1&1&-1&0&-1&1&1&-1 \\ 0&1&-1&0&-1&2&1&-2 \\ 1&1&-1&0&-2&3&1&-2 \\
        1&1&-2&1&-3&4&2&-4 \\ 1&0&-2&1&-2&3&2&-3 \\ 1&0&-2&1&-2&3&1&-2 \\
        0&0&-1&1&-2&2&1&-1 \\ 0&0&-1&1&-1&1&0&0 \end{array}\right)$,
  $s_{9} := \left(\begin{array}   {cccccccc} -1&-2&0&2&-1&0&1&-1 \\ -1&-3&1&2&-2&1&1&-1 \\
        -1&-4&0&4&-3&1&1&-1 \\ -2&-5&1&5&-5&2&2&-2 \\
        -1&-4&0&4&-4&2&2&-2 \\ -1&-3&0&3&-3&1&2&-1 \\
        -1&-2&0&2&-2&1&1&-1 \\ 0&-1&0&1&-1&0&1&-1 \end{array}\right)$,  \\
  $s_{10} := \left(\begin{array}   {cccccccc} -1&1&0&1&-2&1&-1&1 \\ -1&1&-1&2&-2&0&-1&2 \\ -2&2&0&2&-3&0&-1&2 \\
        -3&2&0&3&-4&0&-2&4 \\ -3&2&0&2&-3&0&-1&3 \\ -2&2&0&1&-2&0&-1&2 \\
        -1&2&0&0&-1&0&-1&2 \\ -1&1&0&0&0&0&-1&1 \end{array}\right)$,
  $s_{11} := \left(\begin{array}   {cccccccc} -1&-2&1&1&-1&1&0&0 \\ -3&-2&2&1&-1&1&0&0 \\ -3&-3&2&2&-2&2&-1&1 \\
        -5&-4&3&3&-3&2&0&1 \\ -4&-3&2&2&-2&2&0&1 \\ -3&-2&2&1&-2&2&0&1 \\
        -2&-2&1&1&-1&1&0&1 \\ -1&-1&1&0&0&0&0&1 \end{array}\right)$,
         \\
  $s_{12} := \left(\begin{array}   {cccccccc} 1&0&-1&1&0&-1&1&-1 \\ 1&0&-1&1&0&-1&1&0 \\ 1&-1&-1&2&0&-2&2&-1 \\
        1&-1&-2&3&0&-2&2&-1 \\ 1&-1&-2&2&1&-2&2&-1 \\
        1&-1&-2&2&0&-1&2&-1 \\ 0&0&-1&1&0&-1&2&-1 \\ 0&0&-1&1&0&-1&1&0 \end{array}\right)$,
  $s_{13} := \left(\begin{array}   {cccccccc} 0&0&2&-1&0&-1&0&1 \\ 0&1&3&-2&0&-1&0&0 \\ 0&0&4&-2&0&-1&-1&1 \\
        1&1&5&-3&0&-2&-1&1 \\ 1&1&4&-3&0&-1&-1&1 \\ 1&1&3&-3&1&-1&-1&1 \\
        1&0&2&-2&1&-1&0&0 \\ 1&0&1&-1&0&0&0&0 \end{array}\right)$,
        \\
  $s_{14} := \left(\begin{array}   {cccccccc} 1&0&0&0&0&0&0&0 \\ 1&1&0&0&-1&1&0&-1 \\ 1&0&1&0&-1&1&0&-1 \\
        2&0&0&1&-2&2&0&-2 \\ 2&0&0&0&-1&2&0&-2 \\ 1&0&0&0&-1&2&0&-1 \\
        1&0&0&0&-1&1&1&-1 \\ 1&0&0&0&-1&1&0&0 \end{array}\right)$,
  $s_{15} := \left(\begin{array}   {cccccccc} 0&1&0&0&0&0&0&-1 \\ 1&1&-1&1&-1&0&1&-2 \\ 0&1&0&1&-1&0&1&-3 \\
        1&1&-1&2&-2&1&1&-4 \\ 1&1&-1&1&-1&1&1&-4 \\ 0&1&-1&1&-1&1&1&-3 \\
        0&1&0&0&-1&1&1&-2 \\ 0&1&0&0&-1&1&0&-1 \end{array}\right)$,\\
  $s_{16} := \left(\begin{array}   {cccccccc} -1&1&1&-1&-1&1&1&0 \\ -1&2&1&-1&-2&1&2&-1 \\ -2&2&2&-2&-2&2&2&-1 \\
        -2&3&2&-2&-4&3&3&-1 \\ -2&2&2&-1&-4&3&2&-1 \\
        -2&2&2&-1&-3&2&1&0 \\ -1&1&2&-1&-2&1&1&0 \\ 0&1&1&-1&-1&1&0&0 \end{array}\right)$,
  $s_{17} := \left(\begin{array}   {cccccccc} 0&0&0&-1&2&-1&1&-2 \\1&-1&0&-1&2&-1&2&-3 \\
        1&-1&-1&-1&3&-1&2&-4 \\ 2&-1&-1&-2&4&-1&3&-6 \\
        2&-1&-1&-1&3&-1&2&-5 \\ 1&-1&-1&0&2&-1&2&-4 \\
        0&-1&0&0&1&0&1&-3 \\ 0&0&0&0&0&0&1&-2 \end{array}\right)$,
        \\
  $s_{18} := \left(\begin{array}   {cccccccc} 0&-1&0&0&0&0&0&1 \\ 0&-1&1&-1&0&1&-1&1 \\ 1&-1&0&-1&0&1&-1&2 \\
        1&-1&1&-2&0&1&-1&2 \\ 1&-1&1&-1&-1&1&-1&2 \\
        1&-1&1&-1&0&0&-1&2 \\ 1&-1&0&0&0&0&-1&1 \\ 1&-1&0&0&0&0&0&0 \end{array}\right)$,
  $s_{19} := \left(\begin{array}   {cccccccc} 0&0&1&-1&1&0&0&-2 \\ -1&0&2&-2&2&0&0&-2 \\ -1&0&2&-2&2&1&-1&-3 \\
        -1&0&3&-3&3&1&-1&-4 \\ -1&0&2&-2&3&0&-1&-3 \\
        -1&0&1&-1&2&0&0&-3 \\ -1&-1&1&0&1&0&0&-2 \\ -1&0&1&0&0&0&0&-1
        \end{array}\right)$,  \\
  $s_{20} := \left(\begin{array}   {cccccccc} 1&0&-1&-1&2&0&-1&0 \\ 0&0&-1&-1&3&0&-2&0 \\
        1&1&-2&-2&4&0&-2&0 \\ 1&1&-2&-3&6&0&-3&-1 \\
        1&1&-1&-3&5&0&-3&0 \\ 1&0&-1&-2&4&0&-2&0 \\ 1&0&-1&-1&2&1&-2&0 \\
        0&0&0&-1&1&1&-1&0 \end{array}\right)$,
  $s_{21} := \left(\begin{array}   {cccccccc} 1&1&-1&-1&1&0&1&-1 \\ 0&2&-1&-1&1&0&1&-1 \\ 0&2&-1&-2&2&0&2&-2 \\
        0&3&-2&-2&2&1&2&-3 \\ 0&2&-1&-2&2&1&1&-2 \\ 0&1&-1&-1&1&1&1&-1 \\
        0&1&0&-1&0&1&1&-1 \\ 0&1&0&-1&0&1&0&0 \end{array}\right)$,\\
  $s_{22} := \left(\begin{array}   {cccccccc} -1&0&1&0&-1&1&-1&0 \\ -2&0&2&-1&0&1&-2&1 \\ -2&1&2&-1&-1&2&-3&1 \\
        -3&1&4&-2&-1&2&-4&2 \\ -3&1&3&-1&-1&1&-3&2 \\
        -2&1&2&-1&0&0&-2&1 \\ -1&0&1&0&0&0&-2&1 \\ -1&0&1&0&0&0&-1&0 \end{array}\right)$,
  $s_{23} := \left(\begin{array}   {cccccccc} 0&0&-1&1&-1&2&-1&-1 \\ 1&-1&-2&2&-2&3&-1&-1 \\
        1&0&-2&2&-3&4&-1&-2 \\ 1&-1&-3&3&-4&6&-1&-3 \\
        1&-1&-3&3&-3&4&0&-3 \\ 0&-1&-2&3&-3&3&0&-2 \\
        0&0&-2&2&-2&2&0&-1 \\ 0&0&-1&1&-1&1&0&-1
        \end{array}\right)$,  \\
  $s_{24} := \left(\begin{array}   {cccccccc} -1&0&1&0&-1&1&-1&0 \\ 0&-1&1&0&-1&1&-1&0 \\ 0&0&1&0&-2&2&-2&0 \\
        0&0&2&-1&-2&2&-2&0 \\ 0&0&1&0&-2&1&-1&0 \\ 0&0&1&0&-1&0&-1&0 \\
        0&0&0&0&0&0&-1&0 \\ 0&0&0&0&0&0&0&-1 \end{array}\right)$,
  $s_{25} := \left(\begin{array}   {cccccccc} 1&-1&1&-1&1&0&0&-1 \\ 1&-1&1&-1&1&0&1&-3 \\ 2&-2&1&-1&1&1&0&-3 \\
        3&-2&1&-2&2&1&1&-5 \\ 3&-2&0&-1&2&0&1&-4 \\ 2&-2&0&0&1&0&1&-3 \\
        1&-2&0&0&1&0&1&-2 \\ 1&-1&0&0&0&0&1&-1 \end{array}\right)$,
        \\
  $s_{26} := \left(\begin{array}   {cccccccc} 0&-1&1&0&-1&1&-1&1 \\ 0&-1&1&0&-1&1&-1&0 \\ 1&-1&1&0&-2&2&-2&1 \\
        1&-1&2&-1&-2&2&-2&1 \\ 1&-1&1&0&-2&1&-1&1 \\
        1&-1&1&0&-1&0&-1&1 \\ 1&-1&0&0&0&0&-1&1 \\ 1&-1&0&0&0&0&0&0
        \end{array}\right)$,
  $s_{27} := \left(\begin{array}   {cccccccc} 0&1&2&-2&0&0&0&0 \\ 0&1&3&-3&1&-1&0&1 \\ 0&1&4&-4&1&0&-1&1 \\
        1&2&5&-6&2&-1&-1&2 \\ 1&1&4&-4&1&-1&-1&2 \\ 1&1&3&-3&1&-1&-1&1 \\
        1&0&2&-2&1&-1&0&0 \\ 1&0&1&-1&0&0&0&0 \end{array}\right)$,\\
  $s_{28} := \left(\begin{array}   {cccccccc} 1&1&0&-1&0&1&0&-1 \\ 1&1&0&-1&0&1&0&0 \\ 1&1&1&-2&0&2&0&-1 \\
        2&1&0&-2&0&3&0&-1 \\ 2&0&0&-1&0&2&0&-1 \\ 1&0&0&0&-1&2&0&-1 \\
        1&0&0&0&-1&1&1&-1 \\ 1&0&0&0&-1&1&0&0 \end{array}\right)$,
  $s_{29} := \left(\begin{array}   {cccccccc} 0&0&-1&0&1&0&1&-2 \\ 1&-1&-1&0&1&0&2&-3 \\ 0&-1&-1&0&1&1&2&-4 \\
        1&-1&-2&0&1&1&4&-6 \\ 1&-1&-2&1&0&1&3&-5 \\ 0&-1&-1&1&0&0&3&-4 \\
        0&0&0&0&0&0&2&-3 \\ 0&0&0&0&0&0&1&-2 \end{array}\right)$,
        \\
  $s_{30} := \left(\begin{array}   {cccccccc} 0&1&1&-1&-1&2&-1&0 \\ 0&1&1&-1&-1&2&-1&1 \\ 0&1&1&-1&-2&4&-2&1 \\
        1&2&1&-2&-2&5&-3&2 \\ 1&1&1&-1&-2&4&-3&2 \\ 1&1&0&0&-2&3&-2&1 \\
        1&0&0&0&-1&2&-1&0 \\ 1&0&0&0&-1&1&0&0 \end{array}\right)$,
  $s_{31} := \left(\begin{array}   {cccccccc} 1&1&-1&1&-2&1&0&0 \\ 1&1&-1&1&-2&2&-1&0 \\ 1&2&-2&2&-3&2&-1&0 \\
        2&2&-3&3&-4&3&-2&1 \\ 2&2&-3&2&-3&3&-2&1 \\ 1&2&-2&1&-2&2&-1&1 \\
        0&1&-1&1&-2&2&-1&1 \\ 0&1&0&0&-1&1&-1&1 \end{array}\right)$,
         \\
  $s_{32} := \left(\begin{array}   {cccccccc} 0&0&1&0&-1&0&1&0 \\ 1&1&0&0&-1&0&1&0 \\ 1&0&1&0&-1&0&1&0 \\
        2&1&1&0&-2&0&1&1 \\ 2&1&1&-1&-1&0&1&1 \\ 1&1&1&-1&-1&1&0&1 \\
        0&1&1&-1&0&0&0&1 \\ 0&1&1&-1&0&0&0&0 \end{array}\right)$,
  $s_{33} := \left(\begin{array}   {cccccccc} 0&0&-1&1&0&-1&0&0 \\ -1&-1&0&1&0&-1&0&0 \\ -1&-1&-1&2&0&-2&0&1 \\
        -1&-2&-1&3&-1&-2&0&1 \\ -1&-1&-1&2&-1&-1&0&1 \\
        0&-1&-1&1&0&-1&0&1 \\ 0&-1&-1&1&0&0&-1&1 \\ 0&-1&-1&1&0&0&0&0 \end{array}\right)$,\\
  $s_{34} := \left(\begin{array}   {cccccccc} 0&0&1&0&-1&1&0&-1 \\ 0&0&1&0&-1&1&0&0 \\ 0&-1&1&1&-2&2&0&-1 \\
        1&-1&1&1&-2&2&0&-1 \\ 1&-1&1&0&-1&2&0&-1 \\ 1&0&0&0&-1&2&0&-1 \\
        1&0&0&0&-1&1&1&-1 \\ 1&0&0&0&-1&1&0&0 \end{array}\right)$,
  $s_{35} := \left(\begin{array}   {cccccccc} -1&0&0&1&-1&1&0&-1 \\ 0&0&0&1&-2&2&0&-1 \\ 0&0&-1&2&-2&2&0&-2 \\
        0&0&0&2&-3&3&0&-3 \\ 0&0&0&1&-2&3&0&-3 \\ 0&1&0&0&-1&2&0&-2 \\
        0&1&0&0&-1&1&1&-2 \\ 0&1&0&0&-1&1&0&-1 \end{array}\right)$,
        \\
  $s_{36} := \left(\begin{array}   {cccccccc} 0&0&0&0&-1&2&0&-1 \\ 0&1&1&-1&-1&2&0&-1 \\ -1&1&1&-1&-1&3&0&-2 \\
        -1&1&2&-2&-1&4&0&-2 \\ -1&0&2&-1&-1&3&0&-2 \\
        -1&0&1&0&-1&2&0&-1 \\ -1&0&1&0&-1&1&1&-1 \\ -1&0&1&0&-1&1&0&0 \end{array}\right)$,
  $s_{37} := \left(\begin{array}   {cccccccc} 1&1&1&-2&1&0&0&-1 \\ 3&1&0&-2&1&0&0&-1 \\ 3&1&1&-3&1&1&0&-2 \\
        5&2&1&-4&1&1&0&-2 \\ 4&1&1&-3&1&0&1&-2 \\ 3&1&1&-3&1&0&1&-1 \\
        2&0&1&-2&1&0&0&0 \\ 1&0&1&-1&0&0&0&0 \end{array}\right)$,
        \\
  $s_{38} := \left(\begin{array}   {cccccccc} 0&0&0&-1&2&-1&1&-2 \\ -1&1&1&-2&3&-2&1&-2 \\ -1&1&1&-3&4&-2&2&-3 \\
        -1&2&2&-5&6&-3&2&-4 \\ -1&1&2&-4&5&-2&1&-3 \\
        -1&0&2&-3&4&-2&1&-2 \\ -1&0&2&-2&2&-1&1&-2 \\
        -1&0&1&-1&1&0&0&-1 \end{array}\right)$,
  $s_{39} := \left(\begin{array}   {cccccccc} 0&0&-1&2&-1&0&-1&0 \\ 1&-1&-2&3&-1&0&-1&-1 \\
        1&-1&-2&4&-2&0&-2&0 \\ 1&-2&-3&6&-2&-1&-2&-1 \\
        1&-1&-3&5&-2&-1&-1&-1 \\ 0&-1&-2&4&-2&0&-1&-1 \\
        0&-1&-2&3&-1&0&-1&0 \\ 0&-1&-1&2&-1&0&0&0 \end{array}\right)$,\\
  $s_{40} := \left(\begin{array}   {cccccccc} 0&-1&-1&2&0&-1&0&-1 \\ -1&-2&0&3&-1&-1&0&-1 \\
        -1&-2&0&3&0&-2&0&-2 \\ -1&-4&0&5&-1&-2&0&-3 \\
        -1&-3&0&4&-1&-1&-1&-2 \\ 0&-2&0&3&-1&-1&-1&-1 \\
        0&-1&0&2&-1&0&-1&-1 \\ 0&0&0&1&-1&0&0&-1 \end{array}\right)$,
  $s_{41} := \left(\begin{array}   {cccccccc} 0&-2&0&0&1&0&1&-1 \\ 0&-2&1&-1&2&-1&2&-2 \\ -1&-3&1&-1&3&-1&2&-2 \\
        0&-4&1&-2&5&-2&3&-3 \\ 0&-3&1&-2&4&-2&3&-2 \\
        0&-2&0&-1&3&-2&3&-2 \\ 0&-2&0&0&2&-2&2&-1 \\ 0&-1&0&0&1&-1&1&0
        \end{array}\right)$,  \\
  $s_{42} := \left(\begin{array}   {cccccccc} 0&0&-1&0&1&-1&0&1 \\ -1&-1&0&0&1&-1&0&1 \\ -1&0&0&-1&2&-2&0&2 \\
        -1&-1&0&-1&2&-2&0&3 \\ -1&-1&0&0&1&-2&0&3 \\ 0&-1&0&0&1&-2&0&2 \\
        0&-1&0&0&1&-1&-1&2 \\ 0&-1&0&0&1&-1&0&1 \end{array}\right)$,
  $s_{43} := \left(\begin{array}   {cccccccc} 1&1&-1&0&0&1&-2&0 \\ 2&1&-3&1&0&1&-2&0 \\ 2&2&-3&1&0&1&-3&0 \\
        3&2&-5&2&0&2&-5&1 \\ 2&2&-4&1&1&1&-4&1 \\ 1&2&-3&1&0&1&-3&1 \\
        1&2&-2&0&0&1&-2&1 \\ 0&1&-1&0&0&1&-1&0 \end{array}\right)$,
        \\
  $s_{44} := \left(\begin{array}   {cccccccc} 0&-2&0&1&0&0&0&-1 \\ 0&-3&0&1&0&1&-1&0 \\ 1&-4&-1&2&0&1&-1&-1 \\
        0&-6&-1&3&0&2&-2&-1 \\ 0&-5&-1&2&1&1&-1&-1 \\
        0&-4&-1&2&0&1&0&-1 \\ 0&-2&-1&1&0&1&0&-1 \\ 0&-1&-1&1&0&0&0&0
        \end{array}\right)$,
  $s_{45} := \left(\begin{array}   {cccccccc} -2&0&1&0&1&-2&1&0 \\ -2&-1&2&0&1&-3&2&0 \\ -3&-1&2&0&2&-4&2&0 \\
        -4&-1&4&-1&3&-6&3&0 \\ -3&0&3&-1&2&-5&3&0 \\ -2&0&2&0&1&-4&2&0 \\
        -2&0&2&0&0&-2&1&0 \\ -1&0&1&0&0&-1&1&-1 \end{array}\right)$,\\
  $s_{46} := \left(\begin{array}   {cccccccc} -1&2&1&-1&-1&0&1&0 \\ 0&2&1&-1&-2&0&2&0 \\ -1&3&2&-1&-3&0&2&0 \\
        -1&4&2&-1&-4&-1&4&0 \\ -1&3&2&-1&-3&-1&3&0 \\
        0&3&1&-1&-2&-1&2&0 \\ 0&2&0&0&-2&0&1&0 \\ 0&1&0&0&-1&0&0&1 \end{array}\right)$,
  $s_{47} := \left(\begin{array}   {cccccccc} -2&-1&1&0&1&-1&1&0 \\ -2&-2&2&0&1&-2&2&0 \\ -3&-3&2&0&2&-2&2&0 \\
        -4&-4&4&-1&3&-3&3&0 \\ -3&-3&3&-1&2&-2&3&0 \\
        -2&-3&2&0&1&-1&2&0 \\ -2&-2&2&0&0&0&1&0 \\ -1&-1&1&0&0&0&0&1
        \end{array}\right)$,  \\
  $s_{48} := \left(\begin{array}   {cccccccc} 0&1&-1&1&-1&1&-1&0 \\ -1&1&0&1&-2&2&-1&0 \\
        -1&2&-1&2&-3&2&-1&0 \\ -2&3&0&2&-4&3&-2&0 \\
        -2&3&0&1&-3&3&-2&0 \\ -1&3&0&0&-2&2&-1&0 \\ 0&2&0&0&-2&2&-1&0 \\
        0&1&0&0&-1&1&0&-1 \end{array}\right)$,
  $s_{49} := \left(\begin{array}   {cccccccc} -1&-1&0&1&0&0&-1&0 \\ 0&-1&0&1&0&0&-2&0 \\ -1&-1&0&1&1&-1&-2&0 \\
        -1&-2&1&1&1&-1&-3&0 \\ -1&-1&1&0&1&0&-3&0 \\
        -1&-1&1&0&1&0&-3&1 \\ -1&-1&1&0&1&0&-2&0 \\ -1&0&1&0&0&0&-1&0
        \end{array}\right)$,  \\
  $s_{50} := \left(\begin{array}   {cccccccc} 0&-1&1&0&0&1&-2&0 \\ 0&-2&0&1&0&2&-3&0 \\ 1&-2&0&1&0&2&-4&0 \\
        0&-3&0&2&0&3&-6&0 \\ 0&-3&0&2&0&2&-4&-1 \\ 0&-2&0&1&1&1&-3&-1 \\
        0&-1&0&1&0&1&-2&-1 \\ 0&0&0&0&0&1&-1&-1 \end{array}\right)$,
  $s_{51} := \left(\begin{array}   {cccccccc} 0&1&-1&1&-1&0&0&1 \\ 0&1&0&1&-2&0&0&2 \\ -1&2&0&1&-2&0&0&2 \\
        0&3&0&1&-3&0&0&3 \\ 0&3&0&0&-2&0&0&3 \\ 0&3&0&0&-2&0&0&2 \\
        0&2&0&0&-1&-1&1&1 \\ 0&1&0&0&0&-1&0&1 \end{array}\right)$,\\
  $s_{52} := \left(\begin{array}   {cccccccc} -1&-2&0&1&1&-1&1&-1 \\ -1&-3&1&0&2&-1&1&-1 \\
        -1&-4&0&1&2&-1&2&-2 \\ -1&-6&1&1&3&-2&3&-3 \\
        -1&-5&1&1&2&-1&2&-3 \\ -1&-3&1&0&2&-1&2&-3 \\
        -1&-2&1&0&1&0&1&-2 \\ -1&-1&1&0&0&0&1&-1 \end{array}\right)$,
  $s_{53} := \left(\begin{array}   {cccccccc} 0&0&-1&1&1&-2&1&-1 \\ 1&1&-1&0&2&-3&1&-1 \\ 1&0&-2&1&3&-4&1&-1 \\
        2&0&-2&1&4&-6&2&-2 \\ 1&0&-1&1&3&-5&1&-1 \\ 1&0&-1&1&2&-4&1&-1 \\
        1&0&-1&1&1&-2&0&-1 \\ 0&0&0&0&1&-1&0&-1 \end{array}\right)$,
         \\
  $s_{54} := \left(\begin{array}   {cccccccc} 0&0&-1&0&1&0&-1&2 \\ -1&0&-1&1&0&0&-1&3 \\ -1&0&-2&1&1&0&-1&3 \\
        -2&1&-2&1&1&0&-2&5 \\ -2&1&-1&0&1&0&-1&4 \\ -2&1&-1&0&1&0&-1&3 \\
        -1&1&-1&0&1&-1&0&2 \\ 0&1&-1&0&0&0&0&1 \end{array}\right)$,
  $s_{55} := \left(\begin{array}   {cccccccc} 0&2&0&-1&-1&1&0&1 \\ 0&3&-1&0&-2&1&0&1 \\ -1&4&0&-1&-2&1&0&2 \\
        -1&5&-1&0&-4&2&0&3 \\ -1&4&-1&0&-3&1&1&2 \\ -1&3&0&0&-3&1&1&1 \\
        0&2&0&0&-2&0&1&1 \\ 0&1&0&0&-1&0&0&1 \end{array}\right)$,
        \\
  $s_{56} := \left(\begin{array}   {cccccccc} 1&1&0&0&-1&0&0&-1 \\ 2&1&-1&1&-2&1&-1&-1 \\ 2&2&-1&1&-2&0&-1&-1 \\
        3&2&-1&2&-4&1&-2&-1 \\ 3&1&-1&2&-3&0&-1&-1 \\
        2&1&0&1&-2&0&-1&-1 \\ 2&1&0&0&-1&0&-1&0 \\ 1&0&0&0&0&0&-1&0
        \end{array}\right)$,
  $s_{57} := \left(\begin{array}   {cccccccc} -1&-1&0&0&1&0&0&1 \\ -2&-1&1&-1&2&-1&1&1 \\
        -2&-2&1&-1&2&0&1&1 \\ -3&-2&1&-2&4&-1&2&1 \\
        -3&-1&1&-2&3&0&1&1 \\ -2&-1&0&-1&2&0&1&1 \\ -2&-1&0&0&1&0&1&0 \\
        -1&0&0&0&0&0&1&0 \end{array}\right)$,\\
  $s_{58} := \left(\begin{array}   {cccccccc} -1&-1&2&-1&0&1&0&0 \\ 0&-1&2&-2&0&2&0&0 \\ 0&-2&3&-2&0&2&0&0 \\
        0&-3&4&-3&0&3&1&-1 \\ 0&-2&3&-2&-1&3&1&-1 \\
        0&-2&2&-1&-1&3&0&-1 \\ 0&-1&1&-1&0&2&0&-1 \\ 0&-1&0&0&0&1&0&0 \end{array}\right)$,
  $s_{59} := \left(\begin{array}   {cccccccc} 0&1&-1&0&0&0&1&0 \\ 1&2&-1&-1&0&0&2&-1 \\ 0&3&-1&-1&0&0&2&-1 \\
        1&4&-1&-2&0&0&3&-1 \\ 1&3&-1&-1&0&-1&3&-1 \\
        0&2&0&-1&0&-1&3&-1 \\ 0&2&0&-1&0&-1&2&0 \\ 0&1&0&-1&1&-1&1&0
        \end{array}\right)$,  \\
  $s_{60} := \left(\begin{array}   {cccccccc} -1&0&1&-1&1&1&-2&1 \\ 0&0&0&-1&1&2&-3&2 \\ -1&1&1&-2&1&3&-4&2 \\
        -1&1&1&-3&2&4&-5&3 \\ -1&1&1&-3&2&3&-4&3 \\ -1&0&1&-2&1&3&-3&2 \\
        -1&0&0&-1&1&2&-2&1 \\ 0&0&0&-1&1&1&-1&0 \end{array}\right)$,
  $s_{61} := \left(\begin{array}   {cccccccc} 1&1&-1&1&-1&0&0&-1 \\ 1&0&-1&2&-1&-1&1&-2 \\ 2&0&-2&3&-2&0&0&-2 \\
        2&0&-3&5&-3&-1&1&-3 \\ 2&0&-3&4&-2&-1&1&-2 \\
        1&0&-2&3&-2&0&1&-2 \\ 1&0&-2&2&-1&0&1&-2 \\ 0&0&-1&1&0&0&0&-1
        \end{array}\right)$,  \\
  $s_{62} := \left(\begin{array}   {cccccccc} 1&0&1&-1&0&0&0&-1 \\ 0&1&2&-2&0&0&0&-1 \\ 0&0&3&-2&0&0&0&-2 \\
        0&0&4&-3&0&0&0&-2 \\ 0&0&3&-2&0&-1&1&-2 \\ 0&0&2&-1&-1&0&1&-2 \\
        0&0&1&-1&0&0&1&-2 \\ 0&0&0&0&0&0&0&-1 \end{array}\right)$,
  $s_{63} := \left(\begin{array}   {cccccccc} 0&-1&-1&0&1&1&-1&1 \\ -1&-1&0&-1&2&1&-1&1 \\ 0&-2&-1&-1&3&1&-1&1 \\
        0&-2&-1&-2&4&2&-2&2 \\ 0&-2&0&-2&3&2&-2&2 \\
        0&-1&0&-2&3&1&-2&2 \\ 0&-1&0&-1&2&1&-2&1 \\ 0&0&0&-1&1&1&-1&0 \end{array}\right)$,\\
  $s_{64} := \left(\begin{array}   {cccccccc} -1&0&1&-1&0&1&0&1 \\ 0&0&0&-1&0&2&-1&2 \\ -1&1&1&-2&0&2&-1&3 \\
        -1&1&1&-3&1&3&-2&4 \\ -1&1&1&-3&1&3&-2&3 \\ -1&0&1&-2&1&2&-2&3 \\
        -1&0&1&-1&0&2&-2&2 \\ 0&0&1&-1&0&1&-1&1 \end{array}\right)$,
  $s_{65} := \left(\begin{array}   {cccccccc} 0&-2&0&1&0&0&0&-1 \\ 1&-2&-1&1&1&-1&0&-1 \\ 1&-3&0&1&1&-1&0&-2 \\
        2&-4&-1&1&2&-1&0&-3 \\ 2&-3&-1&1&1&0&-1&-2 \\
        2&-3&-1&1&1&0&-1&-1 \\ 1&-2&-1&1&1&0&-1&-1 \\
        1&-1&-1&1&0&0&0&-1 \end{array}\right)$,  \\
  $s_{66} := \left(\begin{array}   {cccccccc} -1&1&2&-2&0&1&0&0 \\ -2&2&2&-2&-1&2&0&0 \\ -3&2&4&-3&-1&2&0&0 \\
        -4&3&5&-4&-2&4&-1&1 \\ -3&2&4&-3&-2&3&0&1 \\
        -2&1&3&-2&-2&3&0&0 \\ -1&1&2&-2&-1&2&0&0 \\ -1&0&1&-1&0&1&0&0 \end{array}\right)$,
  $s_{67} := \left(\begin{array}   {cccccccc} 0&-1&0&1&-1&1&0&-2 \\ 1&0&-1&1&-1&1&0&-3 \\
        1&-1&0&1&-1&1&0&-4 \\ 2&-1&-1&1&-1&2&0&-6 \\
        2&-1&-1&1&-1&2&-1&-4 \\ 2&-1&-1&1&-1&2&-1&-3 \\
        1&-1&-1&1&0&1&-1&-2 \\ 1&-1&-1&1&0&0&0&-1
        \end{array}\right)$,  \\
  $s_{68} := \left(\begin{array}   {cccccccc} 1&0&0&1&-1&0&-1&0 \\ 2&1&-1&1&-1&0&-1&0 \\ 2&1&0&1&-1&-1&-1&0 \\
        4&1&-1&2&-2&-1&-1&0 \\ 3&0&-1&2&-1&-1&-1&0 \\
        2&0&0&1&-1&0&-1&0 \\ 1&0&0&1&-1&0&-1&1 \\ 1&0&0&0&0&0&-1&1
        \end{array}\right)$,
  $s_{69} := \left(\begin{array}   {cccccccc} -1&0&2&-2&1&0&0&1 \\ -2&0&2&-2&1&0&0&2 \\ -3&0&4&-3&1&0&0&2 \\
        -4&0&5&-4&1&1&-1&4 \\ -3&0&4&-3&0&1&0&3 \\ -2&-1&3&-2&0&1&0&2 \\
        -1&0&2&-2&0&1&0&1 \\ -1&0&1&-1&0&1&0&0 \end{array}\right)$,\\
  $s_{70} := \left(\begin{array}   {cccccccc} 1&0&0&0&0&-1&0&0 \\ 2&-1&-1&1&0&-1&-1&1 \\ 2&0&-1&1&0&-2&-1&1 \\
        3&-1&-1&2&-1&-2&-2&2 \\ 3&-1&-1&2&-1&-2&-1&1 \\
        2&-1&0&1&0&-2&-1&1 \\ 2&0&0&0&0&-1&-1&1 \\ 1&0&0&0&0&0&-1&0 \end{array}\right)$,
  $s_{71} := \left(\begin{array}   {cccccccc} -1&1&0&0&0&0&-1&2 \\ -2&0&1&0&0&0&-1&3 \\ -2&1&0&0&0&0&-1&4 \\
        -3&1&1&0&0&-1&-1&6 \\ -3&1&1&0&0&-1&0&4 \\ -3&1&1&0&0&-1&0&3 \\
        -2&1&1&0&-1&0&0&2 \\ -1&1&1&-1&0&0&0&1 \end{array}\right)$,
        \\
  $s_{72} := \left(\begin{array}   {cccccccc} 0&-1&1&-1&1&0&1&-2 \\ 1&-1&0&-1&1&1&1&-3 \\ 0&-1&1&-2&2&1&1&-4 \\
        1&-2&1&-2&2&2&1&-6 \\ 1&-2&1&-2&2&2&0&-4 \\ 1&-1&1&-2&1&2&0&-3 \\
        1&-1&0&-1&1&1&0&-2 \\ 0&-1&0&0&0&1&0&-1 \end{array}\right)$,
  $s_{73} := \left(\begin{array}   {cccccccc} 1&0&-2&1&0&0&1&-2 \\ 2&1&-3&0&1&0&1&-2 \\ 2&1&-4&1&1&-1&2&-3 \\
        4&1&-6&1&2&-1&2&-4 \\ 3&1&-5&1&2&-1&1&-3 \\ 2&1&-3&0&2&-1&1&-3 \\
        1&1&-2&0&1&0&0&-2 \\ 0&0&-1&0&1&0&0&-1 \end{array}\right)$,
        \\
  $s_{74} := \left(\begin{array}   {cccccccc} 0&-1&0&1&-1&0&2&-2 \\ 1&-1&-1&1&0&-1&3&-3 \\ 0&-2&0&1&0&-1&4&-4 \\
        1&-3&-1&2&0&-2&6&-6 \\ 1&-3&-1&2&0&-1&4&-5 \\
        1&-2&-1&1&0&0&3&-4 \\ 0&-1&-1&1&0&0&2&-3 \\ 0&-1&-1&1&0&0&1&-1
        \end{array}\right)$,
  $s_{75} := \left(\begin{array}   {cccccccc} 0&2&0&-1&0&-1&1&0 \\ 0&3&1&-2&0&-1&1&-1 \\ 0&4&1&-2&-1&-1&1&0 \\
        -1&6&2&-3&-1&-2&2&-1 \\ -1&4&2&-2&-1&-2&2&-1 \\
        0&3&1&-1&-1&-2&2&-1 \\ 0&2&1&-1&0&-2&1&0 \\ 0&1&0&0&0&-1&0&0 \end{array}\right)$,\\
  $s_{76} := \left(\begin{array}   {cccccccc} -1&0&0&-1&2&0&-1&1 \\ -2&0&1&-2&3&-1&0&1 \\
        -2&1&1&-3&4&-1&-1&2 \\ -3&1&1&-4&6&-2&0&2 \\
        -3&1&1&-3&4&-1&0&2 \\ -2&1&0&-2&3&-1&0&2 \\ -2&1&0&-1&2&-1&0&1 \\
        -1&1&0&-1&1&0&0&0 \end{array}\right)$,
  $s_{77} := \left(\begin{array}   {cccccccc} -1&0&0&-1&1&0&1&0 \\ -1&0&0&-1&1&-1&3&-1 \\ -2&1&0&-2&2&-1&3&-1 \\
        -2&1&-1&-2&2&-1&5&-2 \\ -2&1&-1&-1&1&-1&4&-1 \\
        -1&1&-1&-1&1&-1&3&-1 \\ -1&0&-1&0&1&-1&2&-1 \\
        -1&0&0&0&0&0&1&-1 \end{array}\right)$,  \\
  $s_{78} := \left(\begin{array}   {cccccccc} 0&0&-1&2&-1&-1&0&0 \\ 0&-1&-2&3&0&-2&0&0 \\ 0&-1&-2&4&-1&-3&1&0 \\
        1&-2&-4&6&-1&-4&1&0 \\ 1&-2&-3&5&-1&-3&1&-1 \\
        1&-2&-2&4&-1&-2&0&0 \\ 1&-1&-2&3&-1&-1&0&0 \\
        1&0&-1&1&0&-1&0&0 \end{array}\right)$,
  $s_{79} := \left(\begin{array}   {cccccccc} 1&-1&0&0&1&-1&-1&2 \\ 1&-1&0&0&1&-1&-2&3 \\ 2&-2&0&0&1&-1&-2&4 \\
        2&-2&1&-1&2&-2&-3&6 \\ 1&-2&1&-1&2&-1&-3&5 \\
        0&-2&1&0&1&-1&-2&4 \\ 0&-1&1&0&0&0&-2&3 \\ 0&0&0&0&0&0&-1&2
        \end{array}\right)$,  \\
  $s_{80} := \left(\begin{array}   {cccccccc} -1&0&0&0&1&0&-2&2 \\ -2&0&1&0&1&0&-3&2 \\ -2&0&1&0&1&0&-4&4 \\
        -4&0&2&0&2&-1&-5&5 \\ -4&0&2&0&1&0&-4&4 \\ -3&0&1&1&0&0&-3&3 \\
        -2&1&1&0&0&0&-2&2 \\ -1&1&0&0&0&0&-1&1 \end{array}\right)$,
  $s_{81} := \left(\begin{array}   {cccccccc} 1&-1&1&-1&1&-1&1&0 \\ 2&-1&1&-2&2&-1&0&1 \\ 2&-2&1&-2&3&-2&1&1 \\
        3&-3&2&-3&4&-2&0&2 \\ 3&-2&1&-2&3&-2&0&2 \\ 2&-1&1&-2&3&-2&0&1 \\
        1&-1&1&-1&2&-2&0&1 \\ 0&0&1&-1&1&-1&0&1 \end{array}\right)$,\\
  $s_{82} := \left(\begin{array}   {cccccccc} 1&-1&-1&2&-1&0&-1&0 \\ 0&-1&-1&3&-1&-1&-1&0 \\
        0&-2&-1&4&-2&0&-2&0 \\ 0&-3&-2&6&-2&-1&-3&1 \\
        0&-2&-2&5&-2&0&-3&1 \\ 0&-2&-1&4&-2&0&-2&1 \\
        0&-1&-1&3&-2&0&-1&1 \\ 0&0&0&1&-1&0&-1&1 \end{array}\right)$,
  $s_{83} := \left(\begin{array}   {cccccccc} 1&1&-1&-1&0&2&-1&0 \\ 0&2&-1&-2&1&2&-1&0 \\ 0&2&-1&-2&0&4&-2&0 \\
        0&4&-1&-4&1&5&-3&0 \\ 0&3&-1&-3&1&4&-3&0 \\ 0&2&0&-3&1&3&-2&0 \\
        0&2&0&-2&0&2&-1&0 \\ 0&1&0&-1&0&1&0&-1 \end{array}\right)$,
        \\
  $s_{84} := \left(\begin{array}   {cccccccc} 0&2&-1&-1&0&1&0&0 \\ 1&2&-1&-2&0&2&0&0 \\ 0&3&-1&-2&-1&3&0&0 \\
        1&4&-2&-3&-1&4&1&-1 \\ 1&3&-2&-2&-1&3&1&-1 \\
        0&2&-1&-2&0&2&1&-1 \\ 0&1&-1&-1&0&2&0&-1 \\ 0&0&-1&0&0&1&0&0 \end{array}\right)$,
  $s_{85} := \left(\begin{array}   {cccccccc} 0&-2&0&1&1&-2&1&0 \\ 0&-3&0&1&2&-2&0&1 \\ 0&-3&0&1&3&-4&1&1 \\
        0&-5&1&1&4&-5&1&2 \\ 0&-4&1&1&3&-4&1&1 \\ -1&-3&1&1&2&-3&1&1 \\
        -1&-2&1&1&1&-2&0&1 \\ 0&-1&0&1&0&-1&0&1 \end{array}\right)$,
         \\
  $s_{86} := \left(\begin{array}   {cccccccc} 1&1&1&-1&0&-1&0&0 \\ 2&1&1&-1&0&-2&1&0 \\ 3&1&1&-1&0&-2&0&0 \\
        4&2&1&-1&0&-4&1&0 \\ 3&1&1&-1&1&-4&1&0 \\ 2&1&1&-1&1&-3&0&1 \\
        2&1&0&-1&1&-2&0&1 \\ 1&0&0&0&0&-1&0&1 \end{array}\right)$,
  $s_{87} := \left(\begin{array}   {cccccccc} -1&-1&-1&1&0&1&0&0 \\ -2&-1&-1&1&0&2&-1&0 \\ -3&-1&-1&1&0&2&0&0 \\
        -4&-2&-1&1&0&4&-1&0 \\ -3&-1&-1&1&-1&4&-1&0 \\
        -2&-1&-1&1&-1&3&0&-1 \\ -2&-1&0&1&-1&2&0&-1 \\
        -1&0&0&0&0&1&0&-1 \end{array}\right)$,\\
  $s_{88} := \left(\begin{array}   {cccccccc} 0&2&0&-1&-1&2&-1&0 \\ 0&3&0&-1&-2&2&0&-1 \\ 0&3&0&-1&-3&4&-1&-1 \\
        0&5&-1&-1&-4&5&-1&-2 \\ 0&4&-1&-1&-3&4&-1&-1 \\
        1&3&-1&-1&-2&3&-1&-1 \\ 1&2&-1&-1&-1&2&0&-1 \\
        0&1&0&-1&0&1&0&-1 \end{array}\right)$,
  $s_{89} := \left(\begin{array}   {cccccccc} 0&1&1&-1&-1&1&0&1 \\ 0&2&0&-1&-1&1&0&2 \\ 0&3&1&-2&-1&1&0&2 \\
        0&4&1&-3&-1&1&1&3 \\ 0&3&1&-2&-1&0&1&3 \\ 1&2&0&-1&-1&0&1&2 \\
        1&1&0&-1&0&0&0&2 \\ 0&0&0&0&0&0&0&1 \end{array}\right)$,  \\
  $s_{90} := \left(\begin{array}   {cccccccc} 0&0&0&1&-1&-1&2&-2 \\ 0&0&1&0&0&-2&3&-3 \\ 0&-1&1&1&-1&-2&4&-4 \\
        0&-1&1&1&0&-4&6&-6 \\ 0&-1&1&1&0&-3&4&-5 \\ 1&-1&0&1&0&-2&3&-4 \\
        0&-1&0&1&0&-1&2&-3 \\ 0&-1&0&1&0&-1&1&-1 \end{array}\right)$,
  $s_{91} := \left(\begin{array}   {cccccccc} 0&0&1&-1&0&-1&2&0 \\ 0&1&1&-1&-1&-1&3&-1 \\ -1&1&2&-2&0&-2&4&-1 \\
        0&2&2&-3&0&-3&6&-2 \\ 0&2&2&-3&0&-2&4&-1 \\ 0&2&1&-2&0&-2&3&0 \\
        0&1&1&-1&0&-2&2&0 \\ 0&1&1&-1&0&-1&1&0 \end{array}\right)$,
        \\
  $s_{92} := \left(\begin{array}   {cccccccc} 0&-2&1&0&0&1&0&-1 \\ 1&-2&0&0&0&2&0&-2 \\ 0&-3&1&0&0&3&-1&-2 \\
        1&-4&1&0&0&4&-1&-3 \\ 0&-3&1&0&0&3&0&-3 \\ 0&-3&0&1&0&2&0&-2 \\
        0&-2&0&0&1&1&0&-1 \\ 0&-1&0&0&1&0&0&-1 \end{array}\right)$,
  $s_{93} := \left(\begin{array}   {cccccccc} -1&-1&-1&2&-1&1&0&-1 \\ -2&-1&-1&2&-1&2&-1&0 \\
        -2&-2&-2&3&-1&2&0&-1 \\ -3&-3&-3&5&-2&3&-1&0 \\
        -3&-2&-2&4&-2&2&0&0 \\ -2&-2&-1&3&-2&2&0&0 \\
        -2&-2&0&2&-1&1&0&0 \\ -1&-1&0&1&0&0&0&0 \end{array}\right)$,\\
  $s_{94} := \left(\begin{array}   {cccccccc} 0&0&-2&2&-1&0&1&-1 \\ 0&-1&-2&3&-1&-1&2&-2 \\
        1&-1&-4&4&-1&-1&2&-2 \\ 1&-1&-5&6&-2&-2&3&-3 \\
        1&-1&-4&5&-2&-1&2&-2 \\ 0&-1&-3&4&-1&-1&1&-1 \\
        0&-1&-2&3&-1&0&0&-1 \\ 0&0&-1&1&0&0&0&-1 \end{array}\right)$,
  $s_{95} := \left(\begin{array}   {cccccccc} -1&1&1&-2&1&0&1&0 \\ 0&2&0&-3&2&0&1&0 \\ -1&2&1&-4&3&0&1&0 \\
        -1&4&1&-6&4&0&1&1 \\ -1&3&1&-5&4&-1&1&1 \\ -1&2&1&-4&3&0&0&1 \\
        0&1&1&-3&2&0&0&1 \\ 0&1&1&-2&1&0&0&0 \end{array}\right)$,
        \\
  $s_{96} := \left(\begin{array}   {cccccccc} 1&-1&0&-1&2&-1&0&0 \\ 1&-1&0&-1&2&-2&1&0 \\ 1&-1&0&-2&4&-3&1&0 \\
        2&-1&0&-3&5&-4&2&-1 \\ 2&-1&0&-2&4&-4&2&-1 \\
        1&-1&0&-1&3&-3&1&-1 \\ 0&-1&1&-1&2&-2&1&-1 \\
        0&-1&0&0&1&-1&1&-1 \end{array}\right)$,
  $s_{97} := \left(\begin{array}   {cccccccc} -1&1&1&-2&1&0&1&0 \\ 0&1&1&-3&1&1&1&0 \\ 0&2&1&-4&2&0&2&0 \\
        0&3&2&-6&2&1&3&-1 \\ 0&2&2&-5&2&0&3&-1 \\ 0&2&1&-4&2&0&2&-1 \\
        0&1&1&-3&2&0&1&-1 \\ 0&0&0&-1&1&0&1&-1 \end{array}\right)$,
        \\
  $s_{98} := \left(\begin{array}   {cccccccc} 1&-1&0&0&1&-2&1&0 \\ 1&-1&0&0&1&-2&0&1 \\ 2&-1&0&0&1&-3&1&0 \\
        2&-1&0&0&2&-5&1&1 \\ 2&0&0&-1&2&-4&1&1 \\ 1&0&0&-1&2&-3&1&0 \\
        1&0&-1&0&1&-2&1&0 \\ 0&0&0&0&0&-1&1&0 \end{array}\right)$,
  $s_{99} := \left(\begin{array}   {cccccccc} 0&0&1&-2&2&0&-1&1 \\ 0&-1&1&-2&3&0&-2&1 \\ 0&0&2&-4&4&0&-2&1 \\
        0&0&2&-5&6&0&-4&2 \\ 0&0&1&-4&5&0&-3&2 \\ 1&0&0&-3&4&0&-2&1 \\
        1&0&0&-2&3&-1&-1&1 \\ 1&0&0&-1&1&0&-1&1 \end{array}\right)$,\\
  $s_{100} := \left(\begin{array}   {cccccccc} -1&1&1&-1&0&1&-2&1 \\ -1&1&1&-1&0&1&-3&3 \\ -1&2&1&-2&1&1&-4&3 \\
        -2&2&1&-2&1&2&-6&5 \\ -1&1&1&-2&1&2&-5&4 \\ 0&1&0&-1&0&2&-4&3 \\
        0&1&0&-1&0&1&-2&2 \\ 0&1&0&-1&0&1&-1&1 \end{array}\right)$,
  $s_{101} := \left(\begin{array}   {cccccccc} 1&1&-1&0&0&1&-2&0 \\ 2&1&-3&1&0&1&-2&0 \\ 2&2&-3&1&0&1&-3&0 \\
        3&2&-5&2&0&2&-5&1 \\ 2&2&-4&2&-1&2&-4&1 \\ 1&2&-3&1&0&1&-3&1 \\
        1&2&-2&0&0&1&-2&1 \\ 0&1&-1&0&0&1&-1&0 \end{array}\right)$,
        \\
  $s_{102} := \left(\begin{array}   {cccccccc} 0&0&1&-2&1&0&1&0 \\ 1&1&0&-3&2&0&1&0 \\ 1&0&1&-4&3&0&1&0 \\
        2&1&1&-6&4&0&2&-1 \\ 2&1&1&-5&3&0&2&-1 \\ 1&1&1&-4&2&1&1&-1 \\
        0&1&1&-3&2&0&1&-1 \\ 0&0&0&-1&1&0&1&-1 \end{array}\right)$,
  $s_{103} := \left(\begin{array}   {cccccccc} 1&0&0&1&-1&0&-1&0 \\ 1&0&0&1&-1&1&-3&1 \\ 2&-1&0&2&-2&1&-3&1 \\
        2&-1&1&2&-2&1&-5&2 \\ 2&-1&1&1&-1&1&-4&1 \\ 1&-1&1&1&-1&1&-3&1 \\
        1&0&1&0&-1&1&-2&1 \\ 1&0&0&0&0&0&-1&1 \end{array}\right)$,
        \\
  $s_{104} := \left(\begin{array}   {cccccccc} 0&1&0&0&-1&0&0&2 \\ 0&3&0&-1&-1&0&0&2 \\ 0&3&1&-1&-2&0&0&3 \\
        -1&5&1&-1&-3&0&0&4 \\ -1&4&1&-1&-2&0&-1&4 \\
        -1&3&1&-1&-2&1&-1&3 \\ -1&2&1&-1&-1&1&-1&2 \\
        -1&1&1&-1&0&0&0&1 \end{array}\right)$,
  $s_{105} := \left(\begin{array}   {cccccccc} 0&0&1&0&-2&1&1&0 \\ -1&1&2&0&-3&1&1&0 \\ -1&0&3&0&-4&2&1&0 \\
        -2&1&4&0&-6&3&1&1 \\ -2&0&3&1&-5&2&1&1 \\ -1&0&2&1&-4&1&1&1 \\
        -1&0&2&0&-2&0&1&1 \\ 0&0&1&0&-1&0&0&1 \end{array}\right)$,\\
  $s_{106} := \left(\begin{array}   {cccccccc} -1&-1&1&0&0&0&-1&2 \\ -2&-1&1&0&1&-1&-1&2 \\ -2&-1&1&0&1&-1&-2&4 \\
        -3&-1&2&-1&2&-2&-2&5 \\ -2&-1&1&0&1&-2&-1&4 \\
        -1&-1&1&0&1&-2&-1&3 \\ 0&-1&0&0&1&-1&-1&2 \\ 0&-1&0&0&1&-1&0&1 \end{array}\right)$,
  $s_{107} := \left(\begin{array}   {cccccccc} 0&1&0&-2&2&0&-1&1 \\ -1&1&0&-2&2&0&-1&2 \\ 0&1&0&-3&3&0&-2&3 \\
        -1&1&0&-4&5&-1&-2&4 \\ -1&0&0&-3&4&0&-2&3 \\
        0&0&-1&-2&3&0&-1&2 \\ 0&0&-1&-1&2&0&-1&1 \\ 0&0&0&-1&1&0&0&0
        \end{array}\right)$,  \\
  $s_{108} := \left(\begin{array}   {cccccccc} 0&1&0&0&-1&-1&2&0 \\ 0&2&0&0&-1&-2&3&-1 \\ -1&2&1&0&-2&-2&4&-1 \\
        -1&3&0&1&-3&-3&6&-2 \\ -1&2&0&1&-2&-3&5&-2 \\
        0&1&0&1&-2&-2&4&-2 \\ 0&0&0&1&-1&-2&3&-1 \\ 0&0&0&1&-1&-1&1&0 \end{array}\right)$,
  $s_{109} := \left(\begin{array}   {cccccccc} 0&1&-1&0&-1&2&-1&0 \\ 0&3&-1&-1&-1&2&-1&0 \\
        0&3&-1&-1&-2&3&-1&0 \\ 1&5&-2&-2&-2&4&-2&0 \\
        1&4&-2&-1&-2&3&-2&1 \\ 1&3&-2&0&-2&2&-1&0 \\
        0&2&-1&0&-1&1&-1&0 \\ 0&1&-1&0&0&0&0&0 \end{array}\right)$,
        \\
  $s_{110} := \left(\begin{array}   {cccccccc} -2&0&1&0&1&-1&0&-1 \\ -3&1&1&-1&2&-1&0&-1 \\
        -4&1&2&-1&2&-1&-1&-1 \\ -6&1&2&-1&3&-1&-1&-2 \\
        -5&1&2&-1&2&-1&0&-2 \\ -4&0&2&-1&2&-1&0&-1 \\
        -2&0&1&-1&2&-1&0&-1 \\ -1&0&0&0&1&-1&0&0 \end{array}\right)$,
  $s_{111} := \left(\begin{array}   {cccccccc} 0&-1&0&1&0&0&-2&1 \\ -1&-2&1&1&1&-1&-2&1 \\ -1&-3&1&2&0&0&-3&1 \\
        -2&-4&2&3&0&-1&-4&2 \\ -1&-4&1&3&0&-1&-3&2 \\
        0&-3&1&2&0&-1&-2&1 \\ 0&-2&0&2&0&-1&-1&0 \\ 0&-1&0&1&0&-1&0&0 \end{array}\right)$,\\
  $s_{112} := \left(\begin{array}   {cccccccc} 1&-1&-1&1&0&0&-1&0 \\ 0&-2&0&1&0&0&-1&0 \\ 1&-2&-1&2&-1&0&-1&0 \\
        1&-4&-1&3&-1&0&-1&-1 \\ 1&-3&-1&2&-1&1&-1&-1 \\
        1&-2&-1&1&0&0&0&-1 \\ 0&-1&-1&1&0&0&0&-1 \\ 0&-1&-1&1&0&0&0&0
        \end{array}\right)$.

The character table of $G^{s_1} =G$:\\\\
\noindent \begin{tabular}{c|cccccccccccccccccccccc}
  & & & & & & & & & &10 & & & & & & & & & &20 & &\\\hline

$\chi _{1}^{(1)}$&1&1&1&1&1&1&1&1&1&1&1&1&1&1&1&1&1&1&1&1&1&1\\

$\chi _{1}^{(2)}$&1&1&1&1&1&1&1&1&1&1&1&1&1&1&-1&-1&-1&-1&1&1&1&1\\

$\chi _{1}^{(3)}$&8&-8&.&-2&2&.&-1&1&-4&4&.&4&-4&.&-2&2&1&-1&1&-1&-3&3\\

$\chi _{1}^{(4)}$&8&-8&.&-2&2&.&-1&1&-4&4&.&4&-4&.&2&-2&-1&1&1&-1&-3&3\\

$\chi _{1}^{(5)}$&28&28&-1&3&3&4&.&.&10&10&-2&8&8&-4&.&.&.&.&2&2&2&2\\

$\chi _{1}^{(6)}$&28&28&-1&3&3&4&.&.&10&10&-2&8&8&-4&.&.&.&.&2&2&2&2\\

$\chi _{1}^{(7)}$&35&35&.&.&.&-5&.&.&5&5&1&7&7&3&-3&-3&.&.&-2&-2&6&6\\
\end{tabular}

\noindent \begin{tabular}{c|cccccccccccccccccccccc}
  & & & & & & & & & &10 & & & & & & & & & &20 & &\\\hline

$\chi _{1}^{({46})}$&840&-840&.&.&.&.&.&.&-60&60&.&20&-20&.&2&-2&-1&1&5&-5&-3&3\\

$\chi _{1}^{({47})}$&840&840&.&.&.&-40&.&.&30&30&.&16&16&8&.&.&.&.&4&4&8&8\\

$\chi _{1}^{({48})}$&840&840&.&.&.&-40&.&.&30&30&.&16&16&8&.&.&.&.&4&4&8&8\\

$\chi _{1}^{({49})}$&972&972&1&-3&-3&36&.&.&.&.&-2&.&.&108&.&.&.&.&.&.&.&.\\

$\chi _{1}^{({50})}$&972&972&1&-3&-3&36&.&.&.&.&-2&.&.&108&.&.&.&.&.&.&.&.\\

$\chi _{1}^{({51})}$&1008&-1008&.&-2&2&.&-1&1&36&-36&.&8&-8&.&.&.&.&.&2&-2&-6&6\\

$\chi _{1}^{({52})}$&1008&-1008&.&-2&2&.&-1&1&36&-36&.&8&-8&.&.&.&.&.&2&-2&-6&6\\
\end{tabular}

\noindent \begin{tabular}{c|cccccccccccccccccccccc}
  & & & & & & & & & &10 & & & & & & & & & &20 & &\\\hline

$\chi _{1}^{({91})}$&3360&-3360&.&.&.&.&.&.&-60&60&.&-16&16&.&.&.&.&.&2&-2&18&
-18\\

$\chi _{1}^{({92})}$&3360&-3360&.&.&.&.&.&.&-60&60&.&-16&16&.&.&.&.&.&2&-2&18&
-18\\

$\chi _{1}^{({93})}$&4096&4096&.&-4&-4&.&-1&-1&64&64&.&.&.&.&.&.&.&.&.&.&.&.\\

$\chi _{1}^{({94})}$&4096&4096&.&-4&-4&.&-1&-1&64&64&.&.&.&.&.&.&.&.&.&.&.&.\\

$\chi _{1}^{({95})}$&4096&-4096&.&-4&4&.&1&-1&64&-64&.&.&.&.&.&.&.&.&.&.&.&.\\

$\chi _{1}^{({96})}$&4096&-4096&.&-4&4&.&1&-1&64&-64&.&.&.&.&.&.&.&.&.&.&.&.\\

$\chi _{1}^{({97})}$&4200&4200&.&.&.&40&.&.&-30&-30&.&-8&-8&-24&.&.&.&.&-2&-2&-6&
-6\\

$\chi _{1}^{({98})}$&4200&4200&.&.&.&40&.&.&-30&-30&.&-8&-8&-24&.&.&.&.&-2&-2&-6&
-6\\

$\chi _{1}^{({99})}$&4200&-4200&.&.&.&.&.&.&60&-60&.&4&-4&.&2&-2&-1&1&1&-1&-3&3\\

$\chi _{1}^{({100})}$&4200&-4200&.&.&.&.&.&.&60&-60&.&4&-4&.&-2&2&1&-1&1&-1&-3&3\\

$\chi _{1}^{({101})}$&4200&4200&.&.&.&-40&.&.&15&15&.&24&24&104&.&.&.&.&.&.&8&8\\

$\chi _{1}^{({102})}$&4480&4480&-1&-5&-5&64&1&1&-44&-44&.&.&.&-128&.&.&.&.&.&.&
16&16\\

$\chi _{1}^{({103})}$&4536&-4536&.&6&-6&.&.&.&.&.&.&12&-12&.&-2&2&1&-1&3&-3&-9&9\\

$\chi _{1}^{({104})}$&4536&-4536&.&6&-6&.&.&.&.&.&.&12&-12&.&2&-2&-1&1&3&-3&-9&9\\

$\chi _{1}^{({105})}$&4536&4536&1&1&1&-24&1&1&81&81&.&-24&-24&-72&.&.&.&.&.&
.&.&.\\

$\chi _{1}^{({106})}$&5600&-5600&.&.&.&.&.&.&-100&100&.&16&-16&.&.&.&.&.&-2&2&6&
-6\\

$\chi _{1}^{({107})}$&5600&-5600&.&.&.&.&.&.&20&-20&.&-16&16&.&.&.&.&.&-4&4&-12&
12\\

$\chi _{1}^{({108})}$&5600&-5600&.&.&.&.&.&.&20&-20&.&-16&16&.&.&.&.&.&-4&4&-12&
12\\

$\chi _{1}^{({109})}$&5670&5670&1&5&5&6&-1&-1&-81&-81&-2&6&6&-90&.&.&.&.&.&.&.&.\\

$\chi _{1}^{({110})}$&6075&6075&.&.&.&-45&.&.&.&.&1&-9&-9&27&-3&-3&.&.&.&.&.&.\\

$\chi _{1}^{({111})}$&6075&6075&.&.&.&-45&.&.&.&.&1&-9&-9&27&3&3&.&.&.&.&.&.\\

$\chi _{1}^{({112})}$&7168&-7168&.&-2&2&.&-1&1&16&-16&.&.&.&.&.&.&.&.&.&.&.&.\\

\end{tabular}

\noindent \begin{tabular}{c|ccccccccccccccccccccccc}
  & & & & & & & &30 & & & & & & & & & &40 & & & & &\\\hline

$\chi _{1}^{(1)}$&1&1&1&1&1&1&1&1&1&1&1&1&1&1&1&1&1&1&1&1&1&1&1\\

$\chi _{1}^{(2)}$&1&1&1&1&-1&-1&-1&-1&1&1&1&1&-1&-1&-1&-1&1&1&-1&-1&-1&-1&1\\

$\chi _{1}^{(3)}$&5&-5&1&-1&-1&1&6&-6&1&-1&-1&1&1&-1&-2&2&-3&3&3&-3&-1&1&-1\\

$\chi _{1}^{(4)}$&5&-5&1&-1&1&-1&-6&6&1&-1&-1&1&-1&1&2&-2&-3&3&-3&3&1&-1&-1\\

$\chi _{1}^{(5)}$&10&10&.&.&.&.&14&14&1&1&1&1&1&1&-2&-2&5&5&2&2&.&.&-2\\

$\chi _{1}^{(6)}$&10&10&.&.&.&.&-14&-14&1&1&1&1&-1&-1&2&2&5&5&-2&-2&.&.&-2\\

$\chi _{1}^{(7)}$&14&14&.&.&.&.&-21&-21&-1&-1&-1&-1&1&1&-5&-5&3&3&-6&-6&.&.&2\\
\end{tabular}

\noindent \begin{tabular}{c|ccccccccccccccccccccccc}
  & & & & & & & &30 & & & & & & & & & &40 & & & & &\\\hline

$\chi _{1}^{({46})}$&21&-21&.&.&.&.&126&-126&.&.&3&-3&-1&1&-10&10&-3&3&-9&9&-1&1&
-1\\

$\chi _{1}^{({47})}$&-24&-24&.&.&.&.&-84&-84&.&.&3&3&1&1&-20&-20&-1&-1&6&6&-2&-2&
.\\

$\chi _{1}^{({48})}$&-24&-24&.&.&.&.&84&84&.&.&3&3&-1&-1&20&20&-1&-1&-6&-6&2&2&.\\

$\chi _{1}^{({49})}$&.&.&-1&-1&1&1&162&162&.&.&.&.&.&.&18&18&.&.&.&.&.&.&.\\

$\chi _{1}^{({50})}$&.&.&-1&-1&-1&-1&-162&-162&.&.&.&.&.&.&-18&-18&.&.&.&.&.&.&.\\

$\chi _{1}^{({51})}$&90&-90&.&.&.&.&252&-252&.&.&9&-9&-3&3&12&-12&3&-3&.&.&.&.&6\\

$\chi _{1}^{({52})}$&90&-90&.&.&.&.&-252&252&.&.&9&-9&3&-3&-12&12&3&-3&.&.&.&.&6\\
\end{tabular}

\noindent \begin{tabular}{c|ccccccccccccccccccccccc}
  & & & & & & & &30 & & & & & & & & & &40 & & & & &\\\hline

$\chi _{1}^{({91})}$&-6&6&.&.&.&.&-336&336&.&.&12&-12&2&-2&-16&16&.&.&-6&6&-4&4&
2\\

$\chi _{1}^{({92})}$&-6&6&.&.&.&.&336&-336&.&.&12&-12&-2&2&16&-16&.&.&6&-6&4&-4&
2\\

$\chi _{1}^{({93})}$&64&64&1&1&1&1&512&512&1&1&-8&-8&.&.&.&.&.&.&-16&-16&.&.&.\\

$\chi _{1}^{({94})}$&64&64&1&1&-1&-1&-512&-512&1&1&-8&-8&.&.&.&.&.&.&16&16&.&.&.\\

$\chi _{1}^{({95})}$&64&-64&1&-1&1&-1&512&-512&-1&1&-8&8&.&.&.&.&.&.&-16&16&.&.&
.\\

$\chi _{1}^{({96})}$&64&-64&1&-1&-1&1&-512&512&-1&1&-8&8&.&.&.&.&.&.&16&-16&.&.&
.\\

$\chi _{1}^{({97})}$&-30&-30&.&.&.&.&-420&-420&.&.&15&15&-1&-1&-4&-4&3&3&.&.&2&2&
-2\\

$\chi _{1}^{({98})}$&-30&-30&.&.&.&.&420&420&.&.&15&15&1&1&4&4&3&3&.&.&-2&-2&-2\\

$\chi _{1}^{({99})}$&-75&75&.&.&.&.&-210&210&.&.&15&-15&1&-1&-26&26&-3&3&-15&
15&1&-1&7\\

$\chi _{1}^{({100})}$&-75&75&.&.&.&.&210&-210&.&.&15&-15&-1&1&26&-26&-3&3&15&-15&
-1&1&7\\

$\chi _{1}^{({101})}$&-120&-120&.&.&.&.&.&.&.&.&-12&-12&.&.&.&.&-4&-4&.&.&.&.&-8\\

$\chi _{1}^{({102})}$&-80&-80&.&.&.&.&.&.&1&1&-20&-20&.&.&.&.&4&4&.&.&.&.&.\\

$\chi _{1}^{({103})}$&-81&81&.&.&.&.&378&-378&.&.&.&.&.&.&-30&30&.&.&-9&9&-1&1&
-3\\

$\chi _{1}^{({104})}$&-81&81&.&.&.&.&-378&378&.&.&.&.&.&.&30&-30&.&.&9&-9&1&-1&
-3\\

$\chi _{1}^{({105})}$&.&.&.&.&.&.&.&.&.&.&.&.&.&.&.&.&.&.&.&.&.&.&.\\

$\chi _{1}^{({106})}$&-10&10&.&.&.&.&.&.&1&-1&2&-2&.&.&.&.&6&-6&.&.&.&.&2\\

$\chi _{1}^{({107})}$&20&-20&.&.&.&.&-280&280&1&-1&11&-11&-1&1&8&-8&-3&3&20&
-20&.&.&-4\\

$\chi _{1}^{({108})}$&20&-20&.&.&.&.&280&-280&1&-1&11&-11&1&-1&-8&8&-3&3&-20&
20&.&.&-4\\

$\chi _{1}^{({109})}$&.&.&.&.&.&.&.&.&.&.&.&.&.&.&.&.&.&.&.&.&.&.&.\\

$\chi _{1}^{({110})}$&.&.&-1&-1&1&1&-405&-405&.&.&.&.&.&.&27&27&.&.&.&.&.&.&.\\

$\chi _{1}^{({111})}$&.&.&-1&-1&-1&-1&405&405&.&.&.&.&.&.&-27&-27&.&.&.&.&.&.&.\\

$\chi _{1}^{({112})}$&-128&128&.&.&.&.&.&.&-1&1&-32&32&.&.&.&.&.&.&.&.&.&.&.\\
\end{tabular}

\noindent \begin{tabular}{c|ccccccccccccccccccccccccc}
  & & & & & 50& & & & & & & & & & 60& & & & & & & & & &70\\\hline

$\chi _{1}^{(1)}$&1&1&1&1&1&1&1&1&1&1&1&1&1&1&1&1&1&1&1&1&1&1&1&1&1\\

$\chi _{1}^{(2)}$&1&1&1&-1&-1&1&1&1&1&-1&-1&-1&-1&-1&1&1&1&1&1&1&1&1&1&-1&-1\\

$\chi _{1}^{(3)}$&1&4&-4&-4&4&-1&1&-2&2&1&-1&.&2&-2&.&.&.&2&-2&.&.&.&.&.&.\\

$\chi _{1}^{(4)}$&1&4&-4&4&-4&-1&1&-2&2&-1&1&.&-2&2&.&.&.&2&-2&.&.&.&.&.&.\\

$\chi _{1}^{(5)}$&-2&4&4&6&6&.&.&.&.&-2&-2&-2&2&2&-4&.&-1&1&1&2&-1&1&-1&2&-2\\

$\chi _{1}^{(6)}$&-2&4&4&-6&-6&.&.&.&.&2&2&2&-2&-2&-4&.&-1&1&1&2&-1&1&-1&-2&2\\

$\chi _{1}^{(7)}$&2&11&11&-9&-9&.&.&3&3&-2&-2&-1&-1&-1&3&-1&.&2&2&-3&.&-2&.&3&-1\\
\end{tabular}

\noindent \begin{tabular}{c|ccccccccccccccccccccccccc}
  & & & & & 50& & & & & & & & & & 60& & & & & & & & & &70\\\hline

$\chi _{1}^{({46})}$&1&4&-4&-4&4&-1&1&-2&2&5&-5&.&2&-2&.&.&.&-6&6&.&.&.&.&.&.\\

$\chi _{1}^{({47})}$&.&24&24&4&4&.&.&.&.&-2&-2&4&-4&-4&8&.&-1&3&3&-10&-1&-1&-1&4&
4\\

$\chi _{1}^{({48})}$&.&24&24&-4&-4&.&.&.&.&2&2&-4&4&4&8&.&-1&3&3&-10&-1&-1&-1&-4&
-4\\

$\chi _{1}^{({49})}$&.&36&36&-6&-6&.&.&.&.&.&.&18&6&6&12&8&.&.&.&.&.&.&.&6&2\\

$\chi _{1}^{({50})}$&.&36&36&6&6&.&.&.&.&.&.&-18&-6&-6&12&8&.&.&.&.&.&.&.&-6&-2\\

$\chi _{1}^{({51})}$&-6&24&-24&-24&24&2&-2&4&-4&.&.&.&-4&4&.&.&.&.&.&.&.&.&.&.&.\\

$\chi _{1}^{({52})}$&-6&24&-24&24&-24&2&-2&4&-4&.&.&.&4&-4&.&.&.&.&.&.&.&.&.&.&.\\
\end{tabular}

\noindent \begin{tabular}{c|ccccccccccccccccccccccccc}
  & & & & & 50& & & & & & & & & & 60& & & & & & & & & &70\\\hline

$\chi _{1}^{({91})}$&-2&16&-16&-16&16&.&.&.&.&2&-2&.&.&.&.&.&.&-6&6&.&.&.&.&.&.\\

$\chi _{1}^{({92})}$&-2&16&-16&16&-16&.&.&.&.&-2&2&.&.&.&.&.&.&-6&6&.&.&.&.&.&.\\

$\chi _{1}^{({93})}$&.&.&.&.&.&.&.&.&.&.&.&.&.&.&.&.&.&-8&-8&.&.&.&.&.&.\\

$\chi _{1}^{({94})}$&.&.&.&.&.&.&.&.&.&.&.&.&.&.&.&.&.&-8&-8&.&.&.&.&.&.\\

$\chi _{1}^{({95})}$&.&.&.&.&.&.&.&.&.&.&.&.&.&.&.&.&.&-8&8&.&.&.&.&.&.\\

$\chi _{1}^{({96})}$&.&.&.&.&.&.&.&.&.&.&.&.&.&.&.&.&.&-8&8&.&.&.&.&.&.\\

$\chi _{1}^{({97})}$&-2&40&40&20&20&.&.&.&.&-4&-4&4&4&4&8&-8&3&-3&-3&-6&3&1&-1&
-4&4\\

$\chi _{1}^{({98})}$&-2&40&40&-20&-20&.&.&.&.&4&4&-4&-4&-4&8&-8&3&-3&-3&-6&3&1&
-1&4&-4\\

$\chi _{1}^{({99})}$&-7&20&-20&-20&20&-1&1&-2&2&-5&5&.&2&-2&.&.&.&6&-6&.&.&.&.&.&
.\\

$\chi _{1}^{({100})}$&-7&20&-20&20&-20&-1&1&-2&2&5&-5&.&-2&2&.&.&.&6&-6&.&.&.&.
.&.\\

$\chi _{1}^{({101})}$&-8&40&40&.&.&.&.&.&.&.&.&.&.&.&8&-8&2&6&6&-1&2&2&2&.&.\\

$\chi _{1}^{({102})}$&.&.&.&.&.&.&.&.&.&.&.&.&.&.&.&.&-2&4&4&4&4&-2&.&.&.\\

$\chi _{1}^{({103})}$&3&60&-60&20&-20&1&-1&2&-2&-3&3&.&6&-6&.&.&.&.&.&.&.&.&.&.&
.\\

$\chi _{1}^{({104})}$&3&60&-60&-20&20&1&-1&2&-2&3&-3&.&-6&6&.&.&.&.&.&.&.&.&.&.&
.\\

$\chi _{1}^{({105})}$&.&-72&-72&.&.&.&.&.&.&.&.&.&.&.&24&8&.&.&.&9&.&.&.&.&.\\

$\chi _{1}^{({106})}$&-2&-80&80&.&.&2&-2&-8&8&.&.&.&.&.&.&.&.&-4&4&.&.&.&.&.&.\\

$\chi _{1}^{({107})}$&4&-80&80&.&.&.&.&.&.&-4&4&.&8&-8&.&.&.&2&-2&.&.&.&.&.&.\\

$\chi _{1}^{({108})}$&4&-80&80&.&.&.&.&.&.&4&-4&.&-8&8&.&.&.&2&-2&.&.&.&.&.&.\\

$\chi _{1}^{({109})}$&.&-90&-90&.&.&.&.&6&6&.&.&.&.&.&6&6&.&.&.&-9&.&.&.&.&.\\

$\chi _{1}^{({110})}$&.&-45&-45&15&15&.&.&3&3&.&.&-9&-9&-9&-21&-1&.&.&.&.&.&.&.&
3&-1\\

$\chi _{1}^{({111})}$&.&-45&-45&-15&-15&.&.&3&3&.&.&9&9&9&-21&-1&.&.&.&.&.&.&.&
-3&1\\

$\chi _{1}^{({112})}$&.&.&.&.&.&.&.&.&.&.&.&.&.&.&.&.&.&-8&8&.&.&.&.&.&.\\
\end{tabular}

\noindent \begin{tabular}{c|cccccccccccccccccccccccc}
  & & & & & & & & & &80 & & & & & & & & & &90 & & & &\\\hline

$\chi _{1}^{(1)}$&1&1&1&1&1&1&1&1&1&1&1&1&1&1&1&1&1&1&1&1&1&1&1&1\\

$\chi _{1}^{(2)}$&-1&-1&-1&-1&-1&1&1&-1&-1&1&1&-1&1&1&1&1&-1&-1&1&1&1&1&1&1\\

$\chi _{1}^{(3)}$&3&-3&-1&1&.&2&-2&.&.&-2&2&.&.&.&-3&3&1&-1&1&-1&.&.&-2&2\\

$\chi _{1}^{(4)}$&-3&3&1&-1&.&2&-2&.&.&-2&2&.&.&.&-3&3&-1&1&1&-1&.&.&-2&2\\

$\chi _{1}^{(5)}$&5&5&-1&-1&-1&1&1&1&.&2&2&.&.&-2&3&3&1&1&-1&-1&.&.&2&2\\

$\chi _{1}^{(6)}$&-5&-5&1&1&1&1&1&-1&.&2&2&.&.&-2&3&3&-1&-1&-1&-1&.&.&2&2\\

$\chi _{1}^{(7)}$&-3&-3&-1&-1&.&2&2&2&1&1&1&1&-1&1&5&5&1&1&1&1&-1&-1&1&1\\
\end{tabular}

\noindent \begin{tabular}{c|cccccccccccccccccccccccc}
  & & & & & & & & & &80 & & & & & & & & & &90 & & & &\\\hline

$\chi _{1}^{({46})}$&9&-9&-1&1&.&2&-2&.&.&2&-2&.&.&.&5&-5&1&-1&1&-1&1&-1&2&-2\\

$\chi _{1}^{({47})}$&-3&-3&-1&-1&1&3&3&1&.&.&.&.&.&2&-5&-5&-1&-1&-1&-1&1&1&-2&-2\\

$\chi _{1}^{({48})}$&3&3&1&1&-1&3&3&-1&.&.&.&.&.&2&-5&-5&1&1&-1&-1&1&1&-2&-2\\

$\chi _{1}^{({49})}$&.&.&.&.&.&.&.&.&.&2&2&.&.&.&-3&-3&-1&-1&1&1&.&.&.&.\\

$\chi _{1}^{({50})}$&.&.&.&.&.&.&.&.&.&2&2&.&.&.&-3&-3&1&1&1&1&.&.&.&.\\

$\chi _{1}^{({51})}$&-9&9&-1&1&.&.&.&.&.&.&.&.&.&.&-3&3&1&-1&1&-1&.&.&2&-2\\

$\chi _{1}^{({52})}$&9&-9&1&-1&.&.&.&.&.&.&.&.&.&.&-3&3&-1&1&1&-1&.&.&2&-2\\
\end{tabular}

\noindent \begin{tabular}{c|cccccccccccccccccccccccc}
  & & & & & & & & & &80 & & & & & & & & & &90 & & & &\\\hline

$\chi _{1}^{({91})}$&-6&6&.&.&.&2&-2&.&.&.&.&.&.&.&-5&5&-1&1&-1&1&-1&1&2&-2\\

$\chi _{1}^{({92})}$&6&-6&.&.&.&2&-2&.&.&.&.&.&.&.&-5&5&1&-1&-1&1&-1&1&2&-2\\

$\chi _{1}^{({93})}$&8&8&.&.&.&.&.&.&.&.&.&.&.&.&-4&-4&.&.&.&.&-1&-1&.&.\\

$\chi _{1}^{({94})}$&-8&-8&.&.&.&.&.&.&.&.&.&.&.&.&-4&-4&.&.&.&.&-1&-1&.&.\\

$\chi _{1}^{({95})}$&-8&8&.&.&.&.&.&.&.&.&.&.&.&.&4&-4&.&.&.&.&-1&1&.&.\\

$\chi _{1}^{({96})}$&8&-8&.&.&.&.&.&.&.&.&.&.&.&.&4&-4&.&.&.&.&-1&1&.&.\\

$\chi _{1}^{({97})}$&3&3&1&1&-1&1&1&1&.&.&.&.&.&-2&.&.&.&.&.&.&.&.&-2&-2\\

$\chi _{1}^{({98})}$&-3&-3&-1&-1&1&1&1&-1&.&.&.&.&.&-2&.&.&.&.&.&.&.&.&-2&-2\\

$\chi _{1}^{({99})}$&3&-3&-1&1&.&-2&2&.&.&2&-2&.&.&.&.&.&.&.&.&.&.&.&-2&2\\

$\chi _{1}^{({100})}$&-3&3&1&-1&.&-2&2&.&.&2&-2&.&.&.&.&.&.&.&.&.&.&.&-2&2\\

$\chi _{1}^{({101})}$&.&.&.&.&.&-2&-2&.&.&.&.&.&.&-1&.&.&.&.&.&.&.&.&3&3\\

$\chi _{1}^{({102})}$&.&.&.&.&.&.&.&.&.&.&.&.&.&4&.&.&.&.&.&.&.&.&.&.\\

$\chi _{1}^{({103})}$&.&.&.&.&.&.&.&.&.&-2&2&.&.&.&4&-4&.&.&.&.&-1&1&.&.\\

$\chi _{1}^{({104})}$&.&.&.&.&.&.&.&.&.&-2&2&.&.&.&4&-4&.&.&.&.&-1&1&.&.\\

$\chi _{1}^{({105})}$&.&.&.&.&.&.&.&.&.&.&.&.&.&-3&6&6&.&.&-2&-2&.&.&-3&-3\\

$\chi _{1}^{({106})}$&.&.&.&.&.&-4&4&.&.&.&.&.&.&.&.&.&.&.&.&.&.&.&-2&2\\

$\chi _{1}^{({107})}$&1&-1&-1&1&.&2&-2&.&.&.&.&.&.&.&.&.&.&.&.&.&.&.&2&-2\\

$\chi _{1}^{({108})}$&-1&1&1&-1&.&2&-2&.&.&.&.&.&.&.&.&.&.&.&.&.&.&.&2&-2\\

$\chi _{1}^{({109})}$&.&.&.&.&.&.&.&.&.&-2&-2&.&-2&3&.&.&.&.&.&.&.&.&3&3\\

$\chi _{1}^{({110})}$&.&.&.&.&.&.&.&.&1&1&1&-1&3&.&.&.&.&.&.&.&.&.&.&.\\

$\chi _{1}^{({111})}$&.&.&.&.&.&.&.&.&-1&1&1&1&3&.&.&.&.&.&.&.&.&.&.&.\\

$\chi _{1}^{({112})}$&.&.&.&.&.&.&.&.&.&.&.&.&.&.&-8&8&.&.&.&.&2&-2&.&.\\
\end{tabular}

\noindent \begin{tabular}{c|cccccccccccccccccc}
  & & & & & & 100& & & & & & & & & &110 & &\\\hline

$\chi _{1}^{(1)}$&1&1&1&1&1&1&1&1&1&1&1&1&1&1&1&1&1&1\\

$\chi _{1}^{(2)}$&1&1&1&1&-1&-1&-1&-1&-1&-1&-1&-1&-1&-1&-1&-1&1&1\\

$\chi _{1}^{(3)}$&.&.&-2&2&.&.&1&-1&2&-2&-3&3&.&.&-2&2&1&-1\\

$\chi _{1}^{(4)}$&.&.&-2&2&.&.&-1&1&-2&2&3&-3&.&.&2&-2&1&-1\\

$\chi _{1}^{(5)}$&.&.&1&1&-1&-1&-1&-1&2&2&4&4&-1&-1&1&1&-1&-1\\

$\chi _{1}^{(6)}$&.&.&1&1&1&1&1&1&-2&-2&-4&-4&1&1&-1&-1&-1&-1\\

$\chi _{1}^{(7)}$&-1&-1&2&2&.&.&-1&-1&-1&-1&-4&-4&.&.&-2&-2&1&1\\
\end{tabular}

\noindent \begin{tabular}{c|cccccccccccccccccc}
  & & & & & & 100& & & & & & & & & &110 & &\\\hline

$\chi _{1}^{({46})}$&.&.&.&.&.&.&1&-1&-1&1&-3&3&.&.&2&-2&-1&1\\

$\chi _{1}^{({47})}$&.&.&.&.&.&.&1&1&1&1&-2&-2&-3&-3&1&1&1&1\\

$\chi _{1}^{({48})}$&.&.&.&.&.&.&-1&-1&-1&-1&2&2&3&3&-1&-1&1&1\\

$\chi _{1}^{({49})}$&.&.&.&.&.&.&-3&-3&.&.&.&.&.&.&.&.&.&.\\

$\chi _{1}^{({50})}$&.&.&.&.&.&.&3&3&.&.&.&.&.&.&.&.&.&.\\

$\chi _{1}^{({51})}$&.&.&.&.&.&.&-3&3&.&.&.&.&.&.&.&.&-1&1\\

$\chi _{1}^{({52})}$&.&.&.&.&.&.&3&-3&.&.&.&.&.&.&.&.&-1&1\\
\end{tabular}

\noindent \begin{tabular}{c|cccccccccccccccccc}
  & & & & & & 100& & & & & & & & & &110 & &\\\hline

$\chi _{1}^{({91})}$&.&.&.&.&.&.&-1&1&1&-1&.&.&.&.&2&-2&2&-2\\

$\chi _{1}^{({92})}$&.&.&.&.&.&.&1&-1&-1&1&.&.&.&.&-2&2&2&-2\\

$\chi _{1}^{({93})}$&.&.&1&1&-1&-1&2&2&-1&-1&.&.&-4&-4&.&.&.&.\\

$\chi _{1}^{({94})}$&.&.&1&1&1&1&-2&-2&1&1&.&.&4&4&.&.&.&.\\

$\chi _{1}^{({95})}$&.&.&-1&1&-1&1&2&-2&1&-1&.&.&-4&4&.&.&.&.\\

$\chi _{1}^{({96})}$&.&.&-1&1&1&-1&-2&2&-1&1&.&.&4&-4&.&.&.&.\\

$\chi _{1}^{({97})}$&.&.&.&.&.&.&.&.&.&.&-2&-2&-3&-3&-1&-1&1&1\\

$\chi _{1}^{({98})}$&.&.&.&.&.&.&.&.&.&.&2&2&3&3&1&1&1&1\\

$\chi _{1}^{({99})}$&.&.&.&.&.&.&.&.&.&.&3&-3&.&.&-2&2&1&-1\\

$\chi _{1}^{({100})}$&.&.&.&.&.&.&.&.&.&.&-3&3&.&.&2&-2&1&-1\\

$\chi _{1}^{({101})}$&-1&-4&.&.&.&.&.&.&.&.&.&.&.&.&.&.&.&.\\

$\chi _{1}^{({102})}$&.&.&-2&-2&.&.&.&.&.&.&.&.&.&.&.&.&.&.\\

$\chi _{1}^{({103})}$&.&.&.&.&.&.&-2&2&-1&1&-3&3&.&.&.&.&.&.\\

$\chi _{1}^{({104})}$&.&.&.&.&.&.&2&-2&1&-1&3&-3&.&.&.&.&.&.\\

$\chi _{1}^{({105})}$&-1&-4&.&.&.&.&.&.&.&.&.&.&.&.&.&.&.&.\\

$\chi _{1}^{({106})}$&.&.&-2&2&.&.&.&.&.&.&.&.&.&.&.&.&-2&2\\

$\chi _{1}^{({107})}$&.&.&1&-1&-1&1&.&.&.&.&.&.&-4&4&2&-2&-1&1\\

$\chi _{1}^{({108})}$&.&.&1&-1&1&-1&.&.&.&.&.&.&4&-4&-2&2&-1&1\\

$\chi _{1}^{({109})}$&1&-2&.&.&.&.&.&.&.&.&.&.&.&.&.&.&.&.\\

$\chi _{1}^{({110})}$&.&3&.&.&.&.&.&.&.&.&.&.&.&.&.&.&.&.\\

$\chi _{1}^{({111})}$&.&3&.&.&.&.&.&.&.&.&.&.&.&.&.&.&.&.\\

$\chi _{1}^{({112})}$&.&.&2&-2&.&.&.&.&.&.&.&.&.&.&.&.&.&.\\

\end{tabular}


The generators of $G^{s_2}$ are:\\

$ \left(\begin{array}   {cccccccc} 1&1&1&-2&0&0& 1&0 \\
1&1&1&-2&0&-1&2&0 \\ 2&1&2&-3&0&-1&2&0
\\
        2&2&3&-4&-1&-1&3&0 \\ 1&1&3&-3&-1&-1&3&0 \\ 1&1&2&-2&-1&-1&2&1 \\
        1&0&1&-1&-1&0&1&1 \\ 1&0&0&0&-1&0&1&0 \end{array}\right)$.

The representatives of conjugacy classes of   $G^{s_2}$ are:\\

$\left(\begin{array}   {cccccccc} 1&0&0&0&0&0&0&0 \\ 0&1&0&0&
0&0&0&0 \\ 0&0&1&0&0&0&0&0 \\
        0&0&0&1&0&0&0&0 \\ 0&0&0&0&1&0&0&0 \\ 0&0&0&0&0&1&0&0 \\
        0&0&0&0&0&0&1&0 \\ 0&0&0&0&0&0&0&1 \end{array}\right)$,
  $\left(\begin{array}{cccccccc} -1&-1&-1&2&0&0&-1&0 \\ -1&-1&-1&2&0&1&-2&0 \\
        -2&-1&-2&3&0&1&-2&0 \\ -2&-2&-3&4&1&1&-3&0 \\
        -1&-1&-3&3&1&1&-3&0 \\ -1&-1&-2&2&1&1&-2&-1 \\
        -1&0&-1&1&1&0&-1&-1 \\ -1&0&0&0&1&0&-1&0 \end{array}\right)$,
  $\left(\begin{array}{cccccccc} -1&-1&1&0&1&-1&-1&1 \\ 0&-1&1&0&1&-2&-1&2 \\ -1&-1&1&0&2&-3&-1&2 \\
        -1&-1&2&-1&3&-4&-2&4 \\ -1&0&2&-1&2&-3&-2&3 \\
        -1&0&2&-1&2&-3&-1&2 \\ 0&0&1&-1&2&-2&-1&1 \\ 0&0&1&-1&1&-1&0&0 \end{array}\right)$,
  $\left(\begin{array}{cccccccc} -1&-1&2&-1&0&1&0&-1 \\ -1&-1&2&-1&-1&3&-1&-1 \\
        -2&-1&3&-2&0&3&-1&-2 \\ -2&-1&4&-3&-1&5&-1&-3 \\
        -1&-1&3&-2&-1&4&-1&-3 \\ 0&-1&2&-2&0&3&-1&-2 \\
        0&0&2&-2&0&2&-1&-1 \\ 0&0&1&-1&0&1&0&-1 \end{array}\right)$,
  $\left(\begin{array}{cccccccc} -1&0&-1&1&0&0&-1&2 \\ -2&0&0&1&0&0&-2&3 \\ -2&1&-1&1&0&0&-2&4 \\
        -3&1&-1&1&1&0&-4&6 \\ -2&1&-1&1&0&1&-4&5 \\ -2&1&-1&1&0&1&-3&3 \\
        -1&1&-1&1&0&0&-2&2 \\ -1&1&0&0&0&0&-1&1 \end{array}\right)$,
  $\left(\begin{array}{cccccccc} -1&0&0&0&0&0&0&0 \\ 0&-1&0&0&0&0&0&0 \\ 0&0&-1&0&0&0&0&0 \\
        0&0&0&-1&0&0&0&0 \\ 0&0&0&0&-1&0&0&0 \\ 0&0&0&0&0&-1&0&0 \\
        0&0&0&0&0&0&-1&0 \\ 0&0&0&0&0&0&0&-1 \end{array}\right)$,
  $\left(\begin{array}{cccccccc} -1&1&-1&1&0&-1&1&-1 \\ -2&1&0&1&0&-2&2&-1 \\
        -2&2&-1&2&-1&-2&2&-1 \\ -3&2&-1&3&-1&-3&3&-2 \\
        -3&1&-1&3&-1&-2&2&-1 \\ -3&1&0&2&-1&-1&1&-1 \\
        -2&0&0&2&-1&-1&1&-1 \\ -1&0&0&1&0&-1&0&0 \end{array}\right)$,
  $\left(\begin{array}{cccccccc} -1&1&1&0&-1&0&0&-1 \\ -1&2&2&-1&-1&-1&1&-1 \\
        -2&2&2&0&-2&-1&1&-1 \\ -3&3&4&-1&-3&-1&2&-2 \\
        -3&3&3&-1&-2&-1&2&-2 \\ -3&2&3&-1&-1&-1&1&-1 \\
        -2&1&2&-1&0&-1&1&-1 \\ -1&0&1&0&0&-1&1&-1 \end{array}\right)$,
  $\left(\begin{array}{cccccccc} -1&1&2&-2&1&0&-1&0 \\ -2&2&3&-3&1&0&-1&1 \\ -3&2&4&-3&1&0&-2&1 \\
        -4&3&6&-5&1&1&-3&2 \\ -4&2&5&-4&1&1&-2&1 \\ -3&1&4&-3&1&1&-2&1 \\
        -2&0&3&-2&1&0&-1&1 \\ -1&0&2&-1&0&0&0&0 \end{array}\right)$,
  $\left(\begin{array}{cccccccc} 0&-1&0&0&2&-2&1&-1 \\ 1&-2&-1&1&2&-3&2&-1 \\ 1&-2&-1&1&3&-4&2&-2 \\
        2&-3&-2&2&4&-6&3&-2 \\ 1&-3&-1&2&3&-5&2&-1 \\
        1&-2&0&1&2&-4&2&-1 \\ 1&-2&0&1&1&-2&1&-1 \\ 1&-1&0&0&1&-1&0&0 \end{array}\right)$,
  $\left(\begin{array}{cccccccc} 0&-1&0&1&0&0&0&-2 \\ 1&-1&-1&1&0&1&0&-3 \\ 0&-2&-1&2&0&1&0&-4 \\
        1&-3&-2&3&0&1&1&-6 \\ 1&-2&-2&2&1&0&1&-5 \\ 1&-2&-1&1&1&0&1&-4 \\
        0&-1&0&0&1&0&1&-3 \\ 0&-1&0&0&1&0&0&-1 \end{array}\right)$,
  $\left(\begin{array}{cccccccc} 0&-1&1&0&0&-1&2&-1 \\ 0&-1&1&1&-1&-2&3&-1 \\ 0&-1&2&0&0&-3&4&-2 \\
        0&-1&2&1&-1&-4&6&-3 \\ 0&-1&2&1&-1&-3&4&-2 \\
        0&0&2&0&-1&-2&3&-1 \\ 1&0&1&0&-1&-1&2&-1 \\ 1&0&0&0&0&-1&1&0 \end{array}%
\right).$

The character table of $G^{s_2}$:\\
{\small \noindent %
}

\noindent \noindent \noindent where A = E(3)$^2$
  = (-1-ER(-3))/2 = -1-b3,
B = E(5),C = E(15)$^{13}$,D = E(15)$^8$,E = E(5)$^2$,F = E(15),G =
E(15)$^{11}$.

The generators of $G^{s_3}$ are:\\
$\left(\begin{array}{cccccccc} 1&-1&-1&0&1& 0&0&1 \\
1&-2&-2&1&1&1&-1&1 \\ 2&-2&-2&0&2&1&-1&1
\\
        3&-3&-4&1&3&1&-2&2 \\ 3&-3&-3&1&2&1&-2&2 \\ 3&-2&-3&1&1&1&-1&1 \\
        2&-1&-2&1&0&1&-1&1 \\ 1&0&-1&0&0&1&-1&1 \end{array}%
\right)$.

The representatives of conjugacy classes of   $G^{s_3}$ are:\\

$\left(%
\right)$.

The character table of $G^{s_3}$:\\
{\small
\noindent %
 }

\noindent where A = E(3)
  = (-1+ER(-3))/2 = b3,
B = E(5)$^2$,C = E(5)$^4$,D = E(15)$^7$,E = E(15)$^{14}$,F =
E(15)$^{13}$,G = E(15)$^4$.

The generators of $G^{s_4}$ are:\\

$\left(%
\right)$.

The representatives of conjugacy classes of   $G^{s_4}$ are:\\

$\left(%
\right)$.

The character table of $G^{s_4}$:\\
${}_{
 \!\!\!\!\!\!\!_{
\!\!\!\!\!\!\!_{
\!\!\!\!\!\!\!  _{  \small
\noindent %
   }}}}$

\noindent where A = E(5)$^4$,B = E(5)$^3$,C = 2*E(5)$^4$,D =
2*E(5)$^3$,E = 3*E(5)$^4$,F = 3*E(5)$^3$,G = 4*E(5)$^4$,H =
4*E(5)$^3$,I = 5*E(5)$^4$,J = 5*E(5)$^3$,K = 6*E(5)$^4$,L =
6*E(5)$^3$,M = -E(5)-E(5)$^4$
  = (1-ER(5))/2 = -b5,
N = -E(5)-E(5)$^3$,O = -E(5)-E(5)$^2$,P = E(5)+E(5)$^2$+E(5)$^3$,Q =
E(5)+E(5)$^2$+E(5)$^4$.

The generators of $G^{s_5}$ are:\\
$\left(%
\right)$.

The character table of $G^{s_5}$:\\
${}_{
 \!\!\!\!_{
\!\!\!\!\!\!\!_{
\!\!\!\!\!\!\!  _{ \small
   \noindent %
\vskip 0.1cm
 \noindent where A = E(3)$^2$
  = (-1-ER(-3))/2 = -1-b3,
B = 4*E(3)$^2$
  = -2-2*ER(-3) = -2-2i3,
C = 5*E(3)$^2$
  = (-5-5*ER(-3))/2 = -5-5b3,
D = 6*E(3)$^2$
  = -3-3*ER(-3) = -3-3i3,
E = 10*E(3)$^2$
  = -5-5*ER(-3) = -5-5i3,
F = 15*E(3)$^2$
  = (-15-15*ER(-3))/2 = -15-15b3,
G = 20*E(3)$^2$
  = -10-10*ER(-3) = -10-10i3,
H = 24*E(3)$^2$
  = -12-12*ER(-3) = -12-12i3,
I = 30*E(3)$^2$
  = -15-15*ER(-3) = -15-15i3,
J = 36*E(3)$^2$
  = -18-18*ER(-3) = -18-18i3,
K = 40*E(3)$^2$
  = -20-20*ER(-3) = -20-20i3,
L = 45*E(3)$^2$
  = (-45-45*ER(-3))/2 = -45-45b3,
M = 60*E(3)$^2$
  = -30-30*ER(-3) = -30-30i3,
N = 64*E(3)$^2$
  = -32-32*ER(-3) = -32-32i3,
O = 80*E(3)$^2$
  = -40-40*ER(-3) = -40-40i3,
P = 81*E(3)$^2$
  = (-81-81*ER(-3))/2 = -81-81b3,
Q = -2*E(3)$^2$
  = 1+ER(-3) = 1+i3,
R = 3*E(3)$^2$
  = (-3-3*ER(-3))/2 = -3-3b3,
S = E(3)-E(3)$^2$
  = ER(-3) = i3,
T = -2*E(3)-E(3)$^2$
  = (3-ER(-3))/2 = 1-b3,
U = 7*E(3)$^2$
  = (-7-7*ER(-3))/2 = -7-7b3,
V = 8*E(3)$^2$
  = -4-4*ER(-3) = -4-4i3,
W = 9*E(3)$^2$
  = (-9-9*ER(-3))/2 = -9-9b3,
X = 2*E(3)-E(3)$^2$
  = (-1+3*ER(-3))/2 = 1+3b3,
Y = -2*E(3)-3*E(3)$^2$
  = (5+ER(-3))/2 = 3+b3,
Z = E(3)+3*E(3)$^2$
  = -2-ER(-3) = -2-i3,
AA = 5*E(3)+2*E(3)$^2$
  = (-7+3*ER(-3))/2 = -2+3b3,
AB = -2*E(3)+3*E(3)$^2$
  = (-1-5*ER(-3))/2 = -3-5b3,
AC = -5*E(3)-3*E(3)$^2$
  = 4-ER(-3) = 4-i3,
AD = -5*E(3)-8*E(3)$^2$
  = (13+3*ER(-3))/2 = 8+3b3,
AE = E(3)-5*E(3)$^2$
  = 2+3*ER(-3) = 2+3i3,
AF = 5*E(3)+6*E(3)$^2$
  = (-11-ER(-3))/2 = -6-b3,
AG = 5*E(3)-3*E(3)$^2$
  = -1+4*ER(-3) = 3+8b3,
AH = -E(3)-6*E(3)$^2$
  = (7+5*ER(-3))/2 = 6+5b3,
AI = 8*E(3)+3*E(3)$^2$
  = (-11+5*ER(-3))/2 = -3+5b3,
AJ = -3*E(3)+6*E(3)$^2$
  = (-3-9*ER(-3))/2 = -6-9b3,
AK = 3*E(3)+9*E(3)$^2$
  = -6-3*ER(-3) = -6-3i3,
AL = -6*E(3)-9*E(3)$^2$
  = (15+3*ER(-3))/2 = 9+3b3,
AM = 2*E(3)+8*E(3)$^2$
  = -5-3*ER(-3) = -5-3i3,
AN = -8*E(3)-6*E(3)$^2$
  = 7-ER(-3) = 7-i3,
AO = -2*E(3)+6*E(3)$^2$
  = -2-4*ER(-3) = -2-4i3.

The generators of $G^{s_6}$ are:\\

$\left(
\right)$.

The representatives of conjugacy classes of   $G^{s_6}$ are:\\

$\left(%
\right)$.

The character table of $G^{s_6}$:\\
${}_{
 \!\!\!\!\!_{
\!\!\!\!\!\!\!_{
\!\!\!\!\!\!\!  _{ \small
  \noindent %

\noindent where A = E(5)$^4$,B = E(5)$^3$,C = 2*E(5)$^4$,D =
2*E(5)$^3$,E = 3*E(5)$^4$,F = 3*E(5)$^3$ ,G = 4*E(5)$^4$,H =
4*E(5)$^3$,I = 5*E(5)$^4$,J = 5*E(5)$^3$,K = 6*E(5)$^4$,L =
6*E(5)$^3$,M = -E(5)-E(5)$^4$
  = (1-ER(5))/2 = -b5,
N = -E(5)-E(5)$^3$,O = -E(5)-E(5)$^2$ ,P = E(5)+E(5)$^2$+E(5)$^3$ ,Q
= E(5)+E(5)$^2$+E(5)$^4$.

The generators of $G^{s_7}$ are:\\
$\left(%
\right)$.

The representatives of conjugacy classes of   $G^{s_7}$ are:\\

$\left(%
\right)$.

The character table of $G^{s_7}$:\\

\noindent %

The generators of $G^{s_8}$ are:\\

$\left(%
\right)$.

The representatives of conjugacy classes of   $G^{s_8}$ are:\\

$\left(%
\right)$.

The character table of $G^{s_8}$:\\

${}_{
 \!\!\!\!\!_{
\!\!\!\!\!\!\!_{ \!\!\!\!\!\!\!  _{  \small
 \noindent
%

\noindent \noindent where A = E(3)$^2$
  = (-1-ER(-3))/2 = -1-b3,
B = 4*E(3)$^2$
  = -2-2*ER(-3) = -2-2i3,
C = 5*E(3)$^2$
  = (-5-5*ER(-3))/2 = -5-5b3,
D = 6*E(3)$^2$
  = -3-3*ER(-3) = -3-3i3,
E = 10*E(3)$^2$
  = -5-5*ER(-3) = -5-5i3,
F = 15*E(3)$^2$
  = (-15-15*ER(-3))/2 = -15-15b3,
G = 20*E(3)$^2$
  = -10-10*ER(-3) = -10-10i3,
H = 24*E(3)$^2$
  = -12-12*ER(-3) = -12-12i3,
I = 30*E(3)$^2$
  = -15-15*ER(-3) = -15-15i3,
J = 36*E(3)$^2$
  = -18-18*ER(-3) = -18-18i3,
K = 40*E(3)$^2$
  = -20-20*ER(-3) = -20-20i3,
L = 45*E(3)$^2$
  = (-45-45*ER(-3))/2 = -45-45b3,
M = 60*E(3)$^2$
  = -30-30*ER(-3) = -30-30i3,
N = 64*E(3)$^2$
  = -32-32*ER(-3) = -32-32i3,
O = 80*E(3)$^2$
  = -40-40*ER(-3) = -40-40i3,
P = 81*E(3)$^2$
  = (-81-81*ER(-3))/2 = -81-81b3,
Q = 2*E(3)$^2$
  = -1-ER(-3) = -1-i3,
R = 3*E(3)$^2$
  = (-3-3*ER(-3))/2 = -3-3b3,
S = 7*E(3)$^2$
  = (-7-7*ER(-3))/2 = -7-7b3,
T = 8*E(3)$^2$
  = -4-4*ER(-3) = -4-4i3,
U = 9*E(3)$^2$
  = (-9-9*ER(-3))/2 = -9-9b3,
V = E(3)-E(3)$^2$
  = ER(-3) = i3,
W = -E(3)-2*E(3)$^2$
  = (3+ER(-3))/2 = 2+b3,
X = 2*E(3)-E(3)$^2$
  = (-1+3*ER(-3))/2 = 1+3b3,
Y = E(3)+3*E(3)$^2$
  = -2-ER(-3) = -2-i3,
Z = -2*E(3)-3*E(3)$^2$
  = (5+ER(-3))/2 = 3+b3,
AA = 5*E(3)+2*E(3)$^2$
  = (-7+3*ER(-3))/2 = -2+3b3,
AB = -5*E(3)-3*E(3)$^2$
  = 4-ER(-3) = 4-i3,
AC = -2*E(3)+3*E(3)$^2$
  = (-1-5*ER(-3))/2 = -3-5b3,
AD = -5*E(3)-8*E(3)$^2$
  = (13+3*ER(-3))/2 = 8+3b3,
AE = E(3)-5*E(3)$^2$
  = 2+3*ER(-3) = 2+3i3,
AF = 5*E(3)-3*E(3)$^2$
  = -1+4*ER(-3) = 3+8b3,
AG = -E(3)-6*E(3)$^2$
  = (7+5*ER(-3))/2 = 6+5b3,
AH = 8*E(3)+3*E(3)$^2$
  = (-11+5*ER(-3))/2 = -3+5b3,
AI = 5*E(3)+6*E(3)$^2$
  = (-11-ER(-3))/2 = -6-b3,
AJ = -3*E(3)+6*E(3)$^2$
  = (-3-9*ER(-3))/2 = -6-9b3,
AK = 3*E(3)+9*E(3)$^2$
  = -6-3*ER(-3) = -6-3i3,
AL = -6*E(3)-9*E(3)$^2$
  = (15+3*ER(-3))/2 = 9+3b3,
AM = 2*E(3)+8*E(3)$^2$
  = -5-3*ER(-3) = -5-3i3,
AN = -2*E(3)+6*E(3)$^2$
  = -2-4*ER(-3) = -2-4i3,
AO = -8*E(3)-6*E(3)$^2$
  = 7-ER(-3) = 7-i3.

The generators of $G^{s_9}$ are:\\

$\left(
\right)$.

The representatives of conjugacy classes of   $G^{s_9}$ are:\\

$\left(%
\right)$.

The character table of $G^{s_9}$:\\
${}_{
 \!\!\!\!\!_{
\!\!\!\!\!\!\!_{ \!\!\!\!\!\!\!  _{  \small
 \noindent
%

\noindent \noindent where A = -E(5)$^3$,B = -E(5),C = -E(3)
  = (1-ER(-3))/2 = -b3,
D = -E(15)$^{14}$,E = -E(15)$^8$,F = -E(15)$^2$ ,G = -E(15)$^{11}$.

The generators of $G^{s_{10}}$ are:\\

$\left(%
\right)$.

The representatives of conjugacy classes of   $G^{s_{10}}$ are:\\

$\left(%
\right)$.

The character table of $G^{s_{10}}$:\\

${}_{
 \!\!\!\!\!\!\!_{
\!\!\!\!\!\!\!_{ \!\!\!\!\!\!\!  _{  \small
 \noindent
%

\noindent \noindent where A = -E(5)$^3$,B = -E(5),C = -E(3)
  = (1-ER(-3))/2 = -b3,
D = -E(15)$^{14}$,E = -E(15)$^8$,F = -E(15)$^2$,G = -E(15)$^{11}$.

The generators of $G^{s_{11}}$ are:\\

$\left(%
\right)$.

The representatives of conjugacy classes of   $G^{s_{11}}$ are:\\

$\left(%
\right)$.

The character table of $G^{s_{11}}$:\\
${}_{
 \!\!\!\!\!_{
\!\!\!\!\!\!\!_{ \!\!\!\!\!\!\!  _{  \small
 \noindent
%

\noindent \noindent where A = E(5)$^4$,B = E(5)$^3$ ,C =
4*E(5)$^4$,D = 4*E(5)$^3$ ,E = 5*E(5)$^4$,F = 5*E(5)$^3$,G =
6*E(5)$^4$ ,H = 6*E(5)$^3$,I = -2*E(5)$^4$,J = -2*E(5)$^3$.

The generators of $G^{s_{12}}$ are:\\

$\left(%
\right)$.

The representatives of conjugacy classes of   $G^{s_{12}}$ are:\\

$\left(%
\right)$.

The character table of $G^{s_{12}}$:\\
\noindent %

\noindent \noindent where A = E(3)$^2$
  = (-1-ER(-3))/2 = -1-b3,
B = 6*E(3)$^2$
  = -3-3*ER(-3) = -3-3i3,
C = 10*E(3)$^2$
  = -5-5*ER(-3) = -5-5i3,
D = 15*E(3)$^2$
  = (-15-15*ER(-3))/2 = -15-15b3,
E = 20*E(3)$^2$
  = -10-10*ER(-3) = -10-10i3,
F = 24*E(3)$^2$
  = -12-12*ER(-3) = -12-12i3,
G = 30*E(3)$^2$
  = -15-15*ER(-3) = -15-15i3,
H = 60*E(3)$^2$
  = -30-30*ER(-3) = -30-30i3,
I = 64*E(3)$^2$
  = -32-32*ER(-3) = -32-32i3,
J = 80*E(3)$^2$
  = -40-40*ER(-3) = -40-40i3,
K = 81*E(3)$^2$
  = (-81-81*ER(-3))/2 = -81-81b3,
L = 90*E(3)$^2$
  = -45-45*ER(-3) = -45-45i3,
M = -3*E(3)$^2$
  = (3+3*ER(-3))/2 = 3+3b3,
N = 2*E(3)$^2$
  = -1-ER(-3) = -1-i3,
O = -7*E(3)$^2$
  = (7+7*ER(-3))/2 = 7+7b3,
P = -8*E(3)$^2$
  = 4+4*ER(-3) = 4+4i3,
Q = 9*E(3)$^2$
  = (-9-9*ER(-3))/2 = -9-9b3,
R = 4*E(3)
  = -2+2*ER(-3) = 4b3,
S = 12*E(3)
  = -6+6*ER(-3) = 12b3,
T = -16*E(3)
  = 8-8*ER(-3) = 8-8i3,
U = -5*E(3)
  = (5-5*ER(-3))/2 = -5b3.

The generators of $G^{s_{13}}$ are:\\
$\left(
\right)$.

The representatives of conjugacy classes of   $G^{s_{13}}$ are:\\

$\left(%
\right)$.

The character table of $G^{s_{13}}$:\\
${}_{
 \!\!\!\!\!_{
\!\!\!\!\!_{ \!\!\!\!\!  _{  \small
 \noindent
%

\noindent \noindent where A = -E(5)$^3$,B = -E(5),C = 2*E(5)$^3$,D =
2*E(5).

The generators of $G^{s_{14}}$ are:\\
$\left(%
\right)$.

The representatives of conjugacy classes of   $G^{s_{14}}$ are:\\
$\left(%
\right)$.

The character table of $G^{s_{14}}$:\\
\noindent %

The generators of $G^{s_{15}}$ are:\\
$\left(%
\right)$.

The representatives of conjugacy classes of   $G^{s_{15}}$ are:\\

$\left(%
\right)$.

The character table of $G^{s_{15}}$:\\
${}_{
 \!\!\!\!\!_{
\!\!\!\!\!_{ \!\!\!\!\!  _{  \small
 \noindent
%

\noindent \noindent where A = E(3)$^2$
  = (-1-ER(-3))/2 = -1-b3,
B = 5*E(3)$^2$
  = (-5-5*ER(-3))/2 = -5-5b3,
C = 9*E(3)$^2$
  = (-9-9*ER(-3))/2 = -9-9b3,
D = 10*E(3)$^2$
  = -5-5*ER(-3) = -5-5i3,
E = 16*E(3)$^2$
  = -8-8*ER(-3) = -8-8i3,
F = -3*E(3)
  = (3-3*ER(-3))/2 = -3b3,
G = -2*E(3)
  = 1-ER(-3) = 1-i3.

The generators of $G^{s_{16}}$ are:\\
$\left(%
\right)$.

The character table of $G^{s_{16}}$:\\
${}_{
 \!\!\!\!\!_{
\!\!\!\!\!_{ \!\!\!\!\!  _{  \small
 \noindent
%
.

\noindent \noindent where A = -E(5)$^3$,B = -E(5),C = -E(3)
  = (1-ER(-3))/2 = -b3,
D = -E(15)$^{14}$,E = -E(15)$^8$,F = -E(15)$^2$,G = -E(15)$^{11}$.

The generators of $G^{s_{17}}$ are:\\
$\left(%
\right)$.

The character table of $G^{s_{17}}$:\\
${}_{
 \!\!\!\!\!\!\!_{
\!\!\!\!\!\!\!_{ \!\!\!\!\!\!\!  _{  \small
 \noindent
%

\noindent \noindent where A = E(5)$^4$,B = E(5)$^3$,C = 4*E(5)$^4$,D
= 4*E(5)$^3$,E = 5*E(5)$^4$,F = 5*E(5)$^3$,G = 6*E(5)$^4$,H =
6*E(5)$^3$,I = -2*E(5)$^4$,J = -2*E(5)$^3$.

The generators of $G^{s_{18}}$ are:\\
$\left(%
\right)$.

The character table of $G^{s_{18}}$:\\
${}_{
 \!\!\!\!\!\!_{
\!\!\!\!\!\!_{ \!\!\!\!\!\!  _{  \small
 \noindent
%

\noindent \noindent where A = E(3)$^2$
  = (-1-ER(-3))/2 = -1-b3,
B = 6*E(3)$^2$
  = -3-3*ER(-3) = -3-3i3,
C = 10*E(3)$^2$
  = -5-5*ER(-3) = -5-5i3,
D = 15*E(3)$^2$
  = (-15-15*ER(-3))/2 = -15-15b3,
E = 20*E(3)$^2$
  = -10-10*ER(-3) = -10-10i3,
F = 24*E(3)$^2$
  = -12-12*ER(-3) = -12-12i3,
G = 30*E(3)$^2$
  = -15-15*ER(-3) = -15-15i3,
H = 60*E(3)$^2$
  = -30-30*ER(-3) = -30-30i3,
I = 64*E(3)$^2$
  = -32-32*ER(-3) = -32-32i3,
J = 80*E(3)$^2$
  = -40-40*ER(-3) = -40-40i3,
K = 81*E(3)$^2$
  = (-81-81*ER(-3))/2 = -81-81b3,
L = 90*E(3)$^2$
  = -45-45*ER(-3) = -45-45i3,
M = 4*E(3)
  = -2+2*ER(-3) = 4b3,
N = 3*E(3)
  = (-3+3*ER(-3))/2 = 3b3,
O = 2*E(3)
  = -1+ER(-3) = 2b3,
P = -5*E(3)$^2$
  = (5+5*ER(-3))/2 = 5+5b3,
Q = -16*E(3)$^2$
  = 8+8*ER(-3) = 8+8i3,
R = 9*E(3)$^2$
  = (-9-9*ER(-3))/2 = -9-9b3,
S = -7*E(3)$^2$
  = (7+7*ER(-3))/2 = 7+7b3,
T = -8*E(3)$^2$
  = 4+4*ER(-3) = 4+4i3,
U = 12*E(3)$^2$
  = -6-6*ER(-3) = -6-6i3.

The generators of $G^{s_{19}}$ are:\\
$\left(
\right)$.

The representatives of conjugacy classes of   $G^{s_{19}}$ are:\\
$\left(%
\right)$.

The character table of $G^{s_{19}}$:\\
${}_{
 \!\!\!\!\!\!\!_{
\!\!\!\!\!\!\!_{ \!\!\!\!\!\!\!  _{  \small
 \noindent
%

\noindent where A = E(5)$^3$, B = E(5), C = -E(4)
  = -ER(-1) = -i,
D = -E(20)$^{17}$, E = -E(20)$^9$.

The generators of $G^{s_{20}}$ are:\\
$\left(%
\right)$.

The representatives of conjugacy classes of   $G^{s_{20}}$ are:\\
$\left(%
\right)$.

The character table of $G^{s_{20}}$:\\
${}_{
 \!\!\!\!\!\!\!_{
\!\!\!\!\!\!\!_{ \!\!\!\!\!\!\!  _{  \small
 \noindent
%

\noindent \noindent where A = E(5)$^3$,B = E(5),C = 2*E(5)$^3$,D =
2*E(5).

The generators of $G^{s_{21}}$ are:\\
$\left(%
\right)$.

The representatives of conjugacy classes of   $G^{s_{21}}$ are:\\
$\left(%
\right)$.

The character table of $G^{s_{21}}$:\\
 \noindent %

The generators of $G^{s_{22}}$ are:\\
 $\left(%
\right)$.

The character table of $G^{s_{22}}$:\\
\noindent %

\noindent \noindent where A = -E(4)
  = -ER(-1) = -i,
B = -4*E(4)
  = -4*ER(-1) = -4i,
C = -5*E(4)
  = -5*ER(-1) = -5i,
D = -6*E(4)
  = -6*ER(-1) = -6i,
E = -10*E(4)
  = -10*ER(-1) = -10i,
F = -15*E(4)
  = -15*ER(-1) = -15i,
G = -20*E(4)
  = -20*ER(-1) = -20i,
H = -3*E(4)
  = -3*ER(-1) = -3i,
I = 2*E(4)
  = 2*ER(-1) = 2i,
J = -9*E(4)
  = -9*ER(-1) = -9i.

The generators of $G^{s_{23}}$ are:\\
$\left(%
\right)$.

The character table of $G^{s_{23}}$:\\
${}_{
 \!\!\!\!\!_{
\!\!\!\!\!_{ \!\!\!\!\!  _{  \small
 \noindent
%

\noindent \noindent where A = -E(5)$^3$,B = -E(5),C = -2*E(5)$^3$,D
= -2*E(5).

The generators of $G^{s_{24}}$ are:\\
$\left(%
\right)$.

The character table of $G^{s_{24}}$:\\
\noindent %
.

The generators of $G^{s_{25}}$ are:\\
$\left(%
\right)$.

The character table of $G^{s_{25}}$:\\
${}_{
 \!\!\!\!\!_{
\!\!\!\!\!_{ \!\!\!\!\!  _{  \small
 \noindent
%

\noindent \noindent where A = -E(5)$^3$,B = -E(5),C = 2*E(5)$^3$,D =
2*E(5).

The generators of $G^{s_{26}}$ are:\\
$\left(%
\right)$.

The character table of $G^{s_{26}}$:\\
\noindent %

The generators of $G^{s_{27}}$ are:\\
$\left(%
\right)$.

The character table of $G^{s_{27}}$:\\
${}_{
 \!\!\!\!\!\!\!_{
\!\!\!\!\!\!\!_{ \!\!\!\!\!\!\!  _{  \small
 \noindent
%

\noindent \noindent where A = E(5)$^4$,B = E(5)$^3$,C = -E(4)
  = -ER(-1) = -i,
D = E(20)$^9$,E = E(20)$^{17}$.

The generators of $G^{s_{28}}$ are:\\
$\left(%
\right)$.

The character table of $G^{s_{28}}$:\\
\noindent %

\noindent \noindent where A = -E(4)
  = -ER(-1) = -i,
B = -4*E(4)
  = -4*ER(-1) = -4i,
C = -5*E(4)
  = -5*ER(-1) = -5i,
D = -6*E(4)
  = -6*ER(-1) = -6i,
E = -10*E(4)
  = -10*ER(-1) = -10i,
F = -15*E(4)
  = -15*ER(-1) = -15i,
G = -20*E(4)
  = -20*ER(-1) = -20i,
H = 2*E(4)
  = 2*ER(-1) = 2i,
I = -3*E(4)
  = -3*ER(-1) = -3i,
J = 9*E(4)
  = 9*ER(-1) = 9i.

\vskip 0.3cm

\noindent{\large\bf Acknowledgement}:\\ We would like to thank Prof.
N. Andruskiewitsch and Dr. F. Fantino for suggestions and help.

\end{document}